\documentclass[twocolumn,english]{IEEEtran} % XBG: this is for the two-column draft
\usepackage{xtab}
\usepackage{scrextend}
\usepackage{cite}

\usepackage{amsthm}
\usepackage{amsmath,amsfonts, amssymb}

\usepackage{bbm} % for $\mathbbm{1}$ indicator function
\usepackage{bm} % bold fonts for greek letters, better than \bm{}
\usepackage{cite}
\usepackage{graphicx}
\graphicspath{{../fig/}} % all figures go here
\usepackage{epstopdf}

\usepackage{booktabs}
\usepackage{caption}
\usepackage{subcaption}

\usepackage{algorithm,algorithmic}% 
\usepackage{enumerate} % different color text
\usepackage{xcolor}

\usepackage{ragged2e} % \justifying command
\usepackage{url}

\usepackage{multirow}
\usepackage{tikz}
\usetikzlibrary{decorations.pathreplacing,positioning, arrows.meta}

\theoremstyle{plain}

\theoremstyle{definition}

\theoremstyle{remark}

\title{Deferrable Load Scheduling under Demand Charge: A Block Model-Predictive Control Approach}

\author{Lei Yang,~\IEEEmembership{Student~Member,~IEEE},
Xinbo Geng,~\IEEEmembership{Member,~IEEE},
Xiaohong Guan,~\IEEEmembership{Fellow,~IEEE},
Lang~Tong,~\IEEEmembership{Fellow,~IEEE},
\thanks{L.Yang and X. Guan are with the Faculty of Electronic and Information Engineering, Xi'an Jiaotong University. Emails: kobela33@stu.xjtu.edu.cn, xhguan@sei.xjtu.edu.cn. The work of L. Yang and X. Guan is supported by National Key R\&D Program of China (2016YFB0901900). L. Yang, X. Geng, and L. Tong are with the School of Electrical and Computer Engineering, Cornell University. Emails:\{ly392, xg72, lt35\}@cornell.edu. The work of X. Geng and L. Tong is supported in part by the U.S. National Science Foundation under Awards 1816397 and 1809830.} % stops a space
}

\begin{document}
\maketitle

\begin{abstract}
Optimal scheduling of deferrable electrical loads can reshape the aggregated load profile to achieve higher operational efficiency and reliability. This paper studies deferrable load scheduling under demand charge that imposes a penalty on the peak consumption over a billing period.  Such a terminal cost poses challenges in real-time scheduling when demand forecasts are inaccurate. A block model-predictive control approach is proposed by breaking the demand charge into a sequence of stage costs. The problem of charging electric vehicles is used to illustrate the efficacy of the proposed approach. Numerical examples show that the block model-predictive control outperforms the benchmark methods in various settings.

\end{abstract}

\begin{IEEEkeywords}
Demand charge, demand side management, deferrable load scheduling, charging of electric vehicles, model predictive control (MPC).
\end{IEEEkeywords}

% For peer review papers, you can put extra information on the cover
% page as needed:
% \ifCLASSOPTIONpeerreview
% \begin{center} \bfseries EDICS Category: 3-BBND \end{center}
% \fi
%
% For peerreview papers, this IEEEtran command inserts a page break and
% creates the second title. It will be ignored for other modes.
\IEEEpeerreviewmaketitle

%===============================================================%
\section{Introduction} % (fold)
\label{sec:introduction}
%\XG{The block MPC approach actually works for any stochastic control problem with demand charge (a special terminal cost function terminal cost function). }

The scheduling of deferrable loads has shown promising benefits to many distribution system applications, such as flexible load aggregation \cite{chang2013real,hao2014aggregate,chen2014distributional,rahmani2017scheduling}, serving electric vehicle (EV) charging demands \cite{Chen&Tong:12SGC,xu2012scheduling,nguyen2013joint,huang2015matching,Yu&Chen&Tong16CSEE,xu2016dynamic,zhang2017real,le2017plug}, and data center power management \cite{dabbagh2017shaving}. The possibility of ``load shifting" enables the optimal shaping of demand patterns, improving power system flexibility,  and achieving enhanced economic goals.  

For commercial and industrial consumers, a significant challenge in real-time scheduling of deferrable loads is that the demand charge imposed on the maximum consumption over a billing period can be a substantial part of the overall operating cost. Because demand charge is a terminal cost at the end of the billing period, a scheduling decision policy needs to consider the impact of the past and current decisions on the overall peak consumption. Without accurate demand forecasts, deferring current demands to the future may incur substantial demand charge when the deferred demands compete with new demand arrivals for services. The problem becomes even more challenging when demands have completion deadlines with penalties on unmet demands. 
%The scheduling of deferrable load has shown great benefits to many power system applications, such as frequency regulation \cite{hao2014aggregate} and energy balance \cite{nguyen2013joint}.
%The possibility of ``load shifting'' enables shaping demand patterns, that can improve the flexibility of power systems and achieve economic goals. %For example, an operator of an electric vehicle (EV) charging station may switch some EVs' charging process on and off during peak demand periods. 
%For commercial and industrial consumers, a significant challenge in scheduling deferrable load is the demand charge imposed on the maximum consumption over a billing period, which can be a substantial part of the electricity bill. Because demand charge is a terminal cost, a scheduling decision needs to consider costs influenced by future demand. For a deferrable load, a simple heuristic is to flatten the load based on load forecasts. Because the billing period is typically much longer than the scheduling interval, load forecasts into the future are inaccurate. Dynamic pricing schemes such as peak pricing and time-of-use pricing also make load-flattening approaches suboptimal.
 
We study the real-time deferrable load scheduling under demand charge (DLS-DC), where stochastic demands arrive with random energy and service completion requests. Practical examples of DLS-DC include scheduling EV charging in public charging facilities and cloud services in large data centers. In both cases, the demand volume during peak hours could make the demand charge substantial for a service provider.
%However, one major challenge for these accomplishments is the demand charge posed to commercial and industrial consumers (e.g. commercial buildings, energy-intensive corporations). Unlike the (per-unit) electricity price that is charged on energy consumption (in kWh), demand charge is on the highest \emph{average} power consumption (in kW) measured in each consecutive 15-minute interval of a billing period (e.g. one month), which can account for a large part of their monthly electricity bills. Therefore, this pricing scheme would significantly make the scheduling much more difficult due to uncertainties brought by deferrable load during the scheduling process.

%In this paper, we study the problem of deferrable load scheduling with demand charge (DLS-DC), where multiple tasks with stochastic arrival, departure and demands are considered. In particular, we concentrate on maximizing the total expected reward over a finite scheduling horizon.  

\subsection{Related Work}
The scheduling of deferrable loads has drawn much attention over the past few decades. 
%such as charging large numbers of EVs \cite{Chen&Tong:12SGC,xu2012scheduling,huang2015matching,Yu&Chen&Tong16CSEE,yu2018deadline}, controlling residential appliances \cite{chang2013real,rahmani2017scheduling}, and mitigating the variability of renewable supply \cite{subramanian2012optimal,papavasiliou2013large}. 
A popular theme is to formulate the problem as %a (static) mathematical programming in a certainty equivalent setting or a
a dynamic program (DP) for real-time control. To overcome DP's curse of dimensionality, index policies and priority rules have been proposed in \cite{yu2018deadline,xu2016dynamic,j2019}. In particular, for the large scale EV charging with stochastic demands and random arrivals, the Whittle's index policy was developed in \cite{yu2018deadline} and shown to be symptomatically optimal in the light traffic regimes.  In \cite{xu2016dynamic},
%Many of these research concentrated on linear programming (LP) and dynamic programming (DP) based approaches for real-time control strategies. 
%The authors of \cite{gan2013real} proposed a real-time distributed control algorithm to reduce the load variance by shifting the power consumption of deferrable load to periods with high renewable generation. 
%The performance of some simple scheduling policies, such as earliest deadline first (EDF) \cite{liu1973scheduling} and least laxity first (LLF) \cite{mok1983fundamental} has been well studied in the deterministic setting. To take into account the exogenous random process, a few recent papers has developed better index rules through DP. 
Xu et al. established a partial priority rule that prioritized EVs with Less Laxity and Longer remaining Processing time (LLLP). This result was further extended by Jin and Xu in \cite{j2019}, who proposed a complete scheduling rule that prioritized EVs under Less Laxity first with Later Deadline (LLF-LD). %It was proved in \cite{j2019} that LLF-LD was optimal when the cost on non-completion penalty was linear with the unmet demand. It is nontrivial, however, to extend these policies to deal with demand charge. 

Deterministic deadline scheduling policies such as Earliest Deadline First (EDF) \cite{liu1973scheduling} and Least Laxity First (LLF) \cite{mok1983fundamental} have also been considered for deferrable load scheduling. These techniques are optimal regardless of the underlying stochastic models under restricted conditions.  In general, such techniques are suboptimal in either robust or average measures, although finite competitive ratio scheduling exists for some cases \cite{Chen&Tong:12SGC}.

Model-predictive control (MPC) has been widely adopted for online scheduling strategies for its simplicity and favorable performance in many applications.  %Garif et. al. \cite{garifi2018stochastic} proposed a chance constrained MPC algorithm to optimally schedule the controllable appliances while satisfying demand response under uncertainties. Such residential-related context was also studied in \cite{zhang2016model}. 
%For example, the author of \cite{rahmani2017scheduling} proposed a multi-scale MPC to schedule the deferrable appliances that required different decision timescales. 
In the absence of demand charge, MPC has been applied for the real-time scheduling of deferrable loads \cite{rahmani2017scheduling}, including cases involving random arrivals of demands \cite{zhang2017real} and random network topology changes \cite{le2017plug}. MPC can also be utilized for tracking a given pre-scheduled trajectory \cite{di2014electric,chen2014distributional}. Chen et al. \cite{chen2014distributional} showed that the distribution of tracking errors would concentrate around if forecasting errors were bounded.

Demand charge has been in place in the U.S. since the 1900s \cite{neufeld1987price}, and it has been applied to EV charging \cite{zhang2017real}, %\cite{qin2016numerical} 
thermostatic control of commercial buildings \cite{wang2017stochastic,zhang2018optimal}, and data centers \cite{xu2014reducing,dabbagh2017shaving}. 
In \cite{kumar2018stochastic}, Kumar et al. proposed a stochastic MPC approach to schedule stationary batteries that simultaneously served local demands and provided frequency regulation services.  Most relevant to this paper, perhaps, are the work of Jin and Xu \cite{j2019,jin2016optimal} and that by Risbeck and Rawllings in \cite{risbeck2020economic}. In \cite{j2019}, Jin and Xu proposed a priority rule for real-time scheduling of EV charging that had strong performance guarantee, although demand charge was not considered. Demand charge was considered explicitly in \cite{jin2016optimal} for the scheduling of non-deferrable loads by tracking the up-to-date peak power in a DP framework.
%To \cite{} the cost of high demand charge, Jin and Xu \cite{jin2016optimal} proposed to track the up-to-date peak power in a DP model for optimal storage operations with random electricity prices for non-deferrable loads. 
%Most relevant to this paper are method demand charge in MPC-based techniques proposed in \cite{kumar2018stochastic,risbeck2020economic}.  In \cite{kumar2018stochastic}, Kumar et al. proposed a stochastic MPC for battery systems scheduling, where they designed a time-varying parameter to weight demand charge in the objective function.
In \cite{risbeck2020economic}, the authors proposed the so-called economic MPC (EMPC) with a special terminal cost and constraint to track a specific reference trajectory. In the absence of such a reference, the EMPC formulation in \cite{risbeck2020economic} does not apply directly.  %where the authors proposed economic MPC (EMPC) with a constructed terminal cost and constraint  to associate with demand charge. However, their approach depended heavily on the accuracy of the demand forecasts that determined the reference trajectory. 

\subsection{Summary of Results}
Formulating DLS-DC as a stochastic optimal control problem, we propose a block model-predictive control (BMPC) approach by tracking the peak consumption and assessing the impact of the demand charge in each stage. Whereas the idea of tracking the peak consumption was considered in \cite{dabbagh2017shaving,jin2016optimal,risbeck2020economic}, BMPC differs from existing techniques in the specific stage costs used in the optimization, and how the peak consumption is tracked over multiple scheduling intervals.

There are apparent similarities between BMPC and EMPC \cite{risbeck2020economic}; both are derived based on the principle of MPC, and both involve some forms of terminal costs. The main difference is that EMPC solves a deterministic optimization problem aimed to track a reference trajectory. Thus the performance of EMPC depends on the quality of such a reference. BMPC, on the other hand, solves a stochastic one with exogenous random parameters, which does not require a reference trajectory but exploits generically short-term forecasts in a rolling-window fashion. The terminal costs used in the two approaches are, therefore, quite different. Another non-trivial difference is that EMPC assumes that the measurement window in which the maximum consumption is measured matches the scheduling interval. Typically in practice, a demand charge is levied on the \emph {maximum average consumption} within several scheduling intervals.

We conduct numerical simulations to compare BMPC with four state-of-the-art algorithms including nominal MPC (NMPC) without demand charge, EMPC\cite{risbeck2020economic}, EDF \cite{liu1973scheduling} and LLF-LD\cite{j2019}. Numerical results demonstrate that BMPC achieves near-optimal performance in most cases. In particular, BMPC can achieve 10\% more total reward than the second best approach when EV charging requests are stochastic. Comparing with EDF and LLF-LD, BMPC can obtain more than 20\% total reward on average.

The remainder of this paper is organized as follows. We formulate the DLS-DC problem as a stochastic optimal control problem in Section \ref{sec:DLS-DC}. In Section \ref{sec:bmpc}, we develope the BMPC algorithm and present justifications on the terminal cost for the demand charge. In Section \ref{sec:application_ev_charging_scheduling}, an application study on EV charging scheduling is presented. Then we demonstrate the numerical results in Section \ref{sec:numerical_results}. Finally, we conclude the paper in Section \ref{sec:conclusion}.

\section{DLS-DC Model} % (fold)
\label{sec:DLS-DC}
%The DLS-DC is first formulated as a finite-horizon stochastic optimal control problem in Section \ref{sub:problem_formulation}. Demand charge is modeled as a structured terminal reward. Main challenges brought by the demand charge are discussed in Section \ref{sub:difficulties_of_dealing_with_demand_charge}.

% \subsection{Scheduling of Stochastic Load with Demand Charge} % (fold)
% \label{sub:scheduling_of_stochastic_load_with_demand_charge}

%\subsection{General Problem Formulation} % (fold)
%\label{sub:problem_formulation}
%We study the scheduling of deferrable loads under a special \emph{demand charge} pricing scheme. 
%This problem is motivated by many real-world applications, such as data center scheduling and EV charging scheduling (see Section \ref{sub:ev_charging}). This section summarizes the common features among these applications and presents a generic problem formulation.
We present in this section a general stochastic control formulation of DLS-DC. Deferrable loads are flexible demands such that their services can be delayed. By shifting part or all of the demands, DLS-DC has strong inter-temporal dependencies. In Section \ref{sec:application_ev_charging_scheduling}, we consider EV charging as a specific form of deferrable loads.

We model the process of scheduling of deferrable loads by a discrete-time dynamic equation  
%The dynamics of deferrable load are represented by the following time-variant system 
\begin{equation}
\label{eqn:system_state}
x_{t+1} = f_t(x_t, u_t, \xi_t),~t \in \mathcal{T}=\{0,1,\cdots,T-1\},
\end{equation}
where the time index $t$ models the decision interval (or stage), $x_t$ the state vector of deferrable loads that includes the amount of unserved demands\footnote{The state vector may also include other attributes of deferrable loads. For the EV charging problem, the state may also includes the deadline for the completion of EV charging.}, $u_t$ the control vector such as the demands to be served in interval $t$, and $\xi_t$ the exogenous random parameter (such as new arrivals of demands) that influences the state of the deferrable loads in the next stage.
%We assume that time is discrete and indexed by $t$. The finite control horizon is denoted by $\mathcal{T}:=\{0,1,\dots,T-1\}$ with a terminal stage $T$. Variables $\mathbf{x}_t \in \mathbb{R}^n$ and $\mathbf{u}_t \in \mathbb{R}^m$ are the state and control vectors at stage $t$. Exogenous random variables are denoted by $\bm{\xi}_t \in \mathbb{R}^p$. We assume that the initial state $\mathbf{x}_0$ and exogenous random variables $\{\bm{\xi}_t \}_{t \in \mathcal{T}}$ are independent.

We consider the problem of optimal DLS-DC by a sequence of control laws $\{\mu_t\}_{t=0}^{T-1}$, where $\mu_t$ maps the state of deferrable loads $x_t$ and the exogenous input $\xi_t$ to the control $u_t = \mu_t(x_t,\xi_t)$ at stage $t$. The objective is to maximize the expected total reward of scheduling under demand charge defined by the following stochastic optimal control problem \eqref{model:osocp}:
%Let $\mu_t(\cdot,\cdot): \mathbb{R}^n \times \mathbb{R}^p \to \mathbb{R}^m$ denote the control law at stage $t$.
%With a little abuse of notations, we use $\mu_t$ to denote $\mu_t(\mathbf{x}_t,\bm{\xi}_t)$ for simplicity. 
%The DLS-DC is formulated as a stochastic optimal control problem \eqref{model:osocp}, which seeks a series of optimal policies $\{\mu_t^*\}_{t \in \mathcal{T}}$ to maximize total expected reward:
\begin{subequations}
\label{model:osocp}
\begin{align}
\max_{\{\mu_t\}_{t \in \mathcal{T}}} ~& \mathbb{E} \left[ \sum_{t\in \mathcal{T}}{\mathcal{G}_t(\tilde{x}_t, \tilde{u}_t,\xi_t)} - \mathcal{C}(\psi) \right ] \label{eqn:sys_obj}\\
\text{s.t.}~& \tilde{x}_{t+1} = f_t(\tilde{x}_t, \tilde{u}_t,\xi_t), \quad t \in \mathcal{T}, \label{eqn:sys_state}  \\
& h_t(\tilde{x}_t,\tilde{u}_t,\xi_t) \leq 0, \quad t \in \mathcal{T},  \label{eqn:sys_cons}\\
& \tilde{u}_t=\mu_t(\tilde{x}_t,\xi_t), \quad t \in \mathcal{T}, \label{eqn:sys_control} \\
&\psi = \max_{t \in \{0,\ell,2\ell,\dots,T-\ell\}}~\left\{ \frac{1}{\ell} \sum_{\tau = t }^{t+\ell-1}{c(\tilde{x}_\tau,\tilde{u}_\tau,\xi_\tau)}  \right \}, \label{eqn:max-power-def} \\
& \text{initializing}~ \tilde{x}_0=x_0, \label{eqn:sys_1}
%& \hspace{3.7cm} \mathcal{T'}:=\{0,\ell,2\ell,\dots,T-\ell\}. \nonumber
\end{align}
\end{subequations}
where $\tilde{x}_t$ and $\tilde{u}_t$ are decision variables corresponding to the state and control at stage $t$ respectively, $\mathcal{G}_t(\cdot)$ the stage reward function, representing the reward of action of serving demands at stage $t$, $h_t(\cdot)$ a set of constraints on the state and control (e.g. maximum power drawn from the grid), $c(\tilde{x}_t,\tilde{u}_t,\xi_t)$ the total demands served in interval $t$, $\psi$ the maximum average demand in non-(overlapping) $\ell$ consecutive scheduling intervals that are referred to as (average-power) measurement window, and $\mathcal{C}(\psi)$ the terminal cost that models the demand charge. Note that the model defined in \eqref{model:osocp} is applicable to a much broader class of scheduling problems beyond DLS-DC.

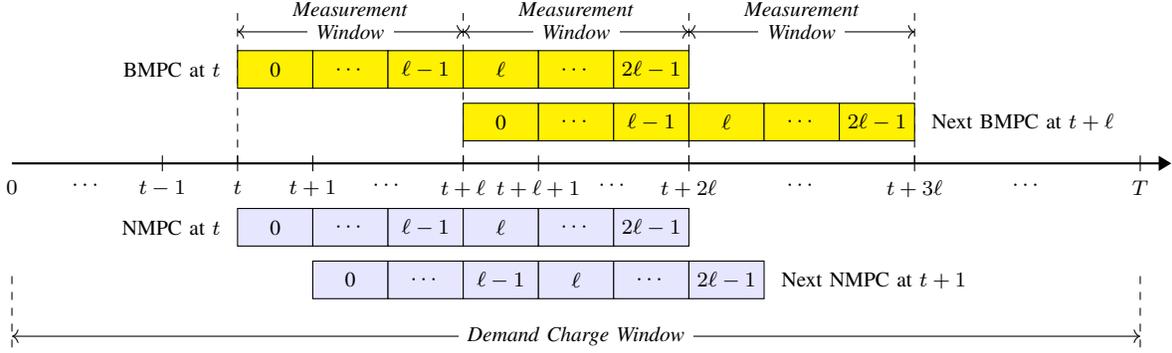
\begin{figure*}[tb]
  \centering
    \begin{tikzpicture}
    \draw[thick, -Triangle] (0,0) -- (0.85*\textwidth,0); %node[font=\footnotesize,right= 3pt]{Horizon};
  
    %draw vertical lines
    \foreach \x in {2,3,4,6,7,9,12,15}
    \draw (\x cm,3pt) -- (\x cm,-3pt);
  
    %draw measurement
    \foreach \x in {3,6,9,12}
    \draw[dashed] (\x cm,0) -- (\x cm,2);
    \draw[<->] (3,1.75) -- (6,1.75) node[midway,font=\footnotesize,fill=white] {\emph{Window}};
    \draw (4.5,2.05) node[font=\footnotesize]{\emph{Measurement}};
    \draw[<->] (6,1.75) -- (9,1.75) node[midway,font=\footnotesize,fill=white] {\emph{Window}};
    \draw (7.5,2.05) node[font=\footnotesize]{\emph{Measurement}};
    \draw[<->] (9,1.75) -- (12,1.75) node[midway,font=\footnotesize,fill=white] {\emph{Window}};
    \draw (10.5,2.05) node[font=\footnotesize] {\emph{Measurement}};
  
    %draw demand charge
    \foreach \x in {0,15}
    \draw[dashed] (\x cm,-1.5) -- (\x cm,-2.5);
    \draw[<->] (0,-2.3) -- (15,-2.3) node[midway,font=\footnotesize,fill=white] {\emph{Demand Charge Window}};
  
    %draw stages
    \draw (0,0) node[below=3pt,font=\footnotesize]{$0$};
    \draw (1,0) node[below=3pt,font=\footnotesize]{$\cdots$};
    \draw (2,0) node[below=3pt,font=\footnotesize]{$t-1$};
    \draw (3,0) node[below=3pt,font=\footnotesize]{$t$};
    \draw (4,0) node[below=3pt,font=\footnotesize]{$t+1$};
    \draw (5,0) node[below=3pt,font=\footnotesize]{$\cdots$};
    \draw (6,0) node[below=3pt,font=\footnotesize]{$t+\ell$};
    \draw (7,0) node[below=3pt,font=\footnotesize]{$t+\ell+1$};
    \draw (8,0) node[below=3pt,font=\footnotesize]{$\cdots$};
  %\draw (9,0) node[below=3pt,font=\footnotesize]{$t+2\ell+1$};
    \draw (9,0) node[below=3pt,font=\footnotesize]{$t+2\ell$};
    \draw (10.5,0) node[below=3pt,font=\footnotesize]{$\cdots$};
    \draw (12,0) node[below=3pt,font=\footnotesize]{$t+3\ell$};
    \draw (13.5,0) node[below=3pt,font=\footnotesize]{$\cdots$};
    \draw (15,0) node[below=3pt,font=\footnotesize]{$T$};
  
    %draw 1st bmpc
    \filldraw[yellow,draw=black] (3,1) rectangle (4,1.5) node[midway,black,font=\footnotesize]{0};
    \filldraw[yellow,draw=black] (4,1) rectangle (5,1.5) node[midway,black,font=\footnotesize]{$\cdots$};
    \filldraw[yellow,draw=black] (5,1) rectangle (6,1.5) node[midway,black,font=\footnotesize]{$\ell-1$};
    \filldraw[yellow,draw=black] (6,1) rectangle (7,1.5)
    node[midway,black,font=\footnotesize]{$\ell$};
    \filldraw[yellow,draw=black] (7,1) rectangle (8,1.5) node[midway,black,font=\footnotesize]{$\cdots$};
    \filldraw[yellow,draw=black] (8,1) rectangle (9,1.5) node[midway,black,font=\footnotesize]{$2\ell-1$};
    \draw (3,1.25) node[left=3pt,font=\footnotesize]{BMPC at $t$};
  
    %draw 1st mpc
    \filldraw[blue!10!white,draw=black] (3,-0.6) rectangle (4,-1.1) node[midway,black,font=\footnotesize]{0};
    \filldraw[blue!10!white,draw=black] (4,-0.6) rectangle (5,-1.1) node[midway,black,font=\footnotesize]{$\cdots$};  \filldraw[blue!10!white,draw=black] (5,-0.6) rectangle (6,-1.1) node[midway,black,font=\footnotesize]{$\ell-1$};
    \filldraw[blue!10!white,draw=black] (6,-0.6) rectangle (7,-1.1) node[midway,black,font=\footnotesize]{$\ell$};
    \filldraw[blue!10!white,draw=black] (7,-0.6) rectangle (8,-1.1) node[midway,black,font=\footnotesize]{$\cdots$};
    \filldraw[blue!10!white,draw=black] (8,-0.6) rectangle (9,-1.1) node[midway,black,font=\footnotesize]{$2\ell-1$};
    \draw (3,-0.85) node[left=3pt,font=\footnotesize]{NMPC at $t$};
  
    %draw 2nd bmpc
    \filldraw[yellow,draw=black] (6,0.3) rectangle (7,0.8) node[midway,black,font=\footnotesize]{0};
    \filldraw[yellow,draw=black] (7,0.3) rectangle (8,0.8) node[midway,black,font=\footnotesize]{$\cdots$};  \filldraw[yellow,draw=black] (8,0.3) rectangle (9,0.8) node[midway,black,font=\footnotesize]{$\ell-1$};
    \filldraw[yellow,draw=black] (9,0.3) rectangle (10,0.8) node[midway,black,font=\footnotesize]{$\ell$};
    \filldraw[yellow,draw=black] (10,0.3) rectangle (11,0.8) node[midway,black,font=\footnotesize]{$\cdots$};
    \filldraw[yellow,draw=black] (11,0.3) rectangle (12,0.8) node[midway,black,font=\footnotesize]{$2\ell-1$};
    \draw (12,0.55) node[right=3pt,font=\footnotesize]{Next BMPC at $t+\ell$};
  
    %draw 2nd mpc
    \filldraw[blue!10!white,draw=black] (4,-1.3) rectangle (5,-1.8) node[midway,black,font=\footnotesize]{0};
    \filldraw[blue!10!white,draw=black] (5,-1.3) rectangle (6,-1.8) node[midway,black,font=\footnotesize]{$\cdots$};
    \filldraw[blue!10!white,draw=black] (6,-1.3) rectangle (7,-1.8) node[midway,black,font=\footnotesize]{$\ell-1$};
    \filldraw[blue!10!white,draw=black] (7,-1.3) rectangle (8,-1.8) node[midway,black,font=\footnotesize]{$\ell$};
    \filldraw[blue!10!white,draw=black] (8,-1.3) rectangle (9,-1.8) node[midway,black,font=\footnotesize]{$\cdots$};
    \filldraw[blue!10!white,draw=black] (9,-1.3) rectangle (10,-1.8) node[midway,black,font=\footnotesize]{$2\ell-1$};
    \draw (10,-1.55) node[right=3pt,font=\footnotesize]{Next NMPC at $t+1$};
    \end{tikzpicture}
  \caption{Temporal structure of BMPC and NMPC with an example of $W=2\ell$.}
  \label{f:bmpc}
\end{figure*}

The main difficulties of dealing with demand charges in the stochastic optimal control framework come from the \emph{mismatch of different timescales}. In particular, three timescales coexist in the formulated DLS-DC model:
\begin{enumerate}
\item the control $u_t=\mu_t(x_t,\xi_t)$ happens at every stage $t \in \mathcal{T} = \{0,1,\cdots,T-1\}$;
\item the peak average consumption that sets the demand charge is calculated at the end of each measurement window at $t = \{\ell,2\ell,\cdots,T\}$;
\item the demand charge is imposed at the terminal stage $T$.
\end{enumerate}
Because a stage reward $\mathcal{G}_t(x_t,u_t,\xi_t)$ is realized at every stage $t$ whereas the demand charge $\mathcal{C}(\psi)$ is levied at the end of the entire control horizon, it is challenging to balance the trade-off between the immediate stage reward and the uncertain demand charge set at the end of the scheduling horizon.  
%It is clearly suboptimal to directly consider $\mathcal{C}(\psi)$ at every stage $t$, due to the fact that $\mathcal{G}_t(x_t,u_t,\xi_t) \ll \mathcal{C}(\psi)$ when demand charge accounts for a significant portion of the bill. Since the control horizon $T$ is typically long, the overall stochastic optimal control problem \eqref{model:osocp} is large-scale and impossible to be directly solved in practice. Hence,
%we have to average or break down demand charge $\mathcal{C}(\psi)$ to every stage $t$. In this case, however, timescale mismatch still exists. The average consumption is evaluated on a \emph{block} of $\ell$ stages, thus a demand charge related cost occurs every $\ell$ stages while the stage reward $\mathcal{G}_t(x_t,u_t,\xi_t)$ happens at every stage $t$.

% subsection difficulties_of_dealing_with_demand_charge (end)

%===============================================================%
\section{Block Model Predictive Control}
\label{sec:bmpc}
%As discussed in Section \ref{sec:DLS-DC},
%the temporal structure of the demand charge pricing scheme suggests some ``block structures'' of control policies. 
This section formalizes a block decision structure and proposes an MPC-based framework for DLS-DC. The proposed framework is termed \emph{Block} MPC because  the rolling window moves a \emph{block} of $\ell$ stages at a time.

% Due to the challenges induced by demand charge, we 
% Generally, the stochastic model formulated in Section \ref{sec:scheduling_of_stochastic_load_with_demand_charge} is not tractable in real-time due to the random system inputs. This future information is not known exactly, but can be forecast so that MPC becomes a natural approach for online solutions. \XG{Better explanation why this problem is hard.} In this section, we formally establish the Block MPC problem with analysis on demand charge.

% \XG{Say sth about the challenges of solving this problem.}

\subsection{Block Model-Predictive Control under Demand Charge} % (fold)
\label{sub:block_mdel_predictive_control}
%\subsubsection{Augmented System States} % (fold)
%\label{ssub:augmented_system_states}
To address the timescale mismatch issues (see Section \ref{sec:DLS-DC}) arising from the demand charge, we introduce an additional system state $\phi_{t}$ at every stage $t$, which tracks the highest average consumption over an $\ell$-sized measurement window \emph{until} stage $t$. The new state variable $\phi_{t}$ evolves according to
% \begin{equation}
% \label{eqn:new_state_update}
% \phi_{t+1} = 
% \begin{cases}
% \frac{1}{\ell}\sum_{\tau=t-\ell+1}^{t}{h(\bm{x}_\tau,\bm{u}_\tau},\bm{\xi}_\tau),&\text{if}~t \in \mathbb{T}; \\
% \phi_{t},&\text{otherwise}.
% \end{cases}
%     % \phi_{t+\ell}=\max\left(\phi_t,\frac{1}{\ell}\sum_{\tau=t}^{t+\ell-1}{h(\bm{x}_\tau,\bm{u}_\tau},\bm{\xi}_\tau)\right), \quad \forall t \in \mathbb{T}
% \end{equation}
\begin{equation}
\label{eqn:new_state_update}
\phi_{t+1} = 
\begin{cases}
\max\left\{\phi_{t},\frac{\sum_{\tau=t-\ell+1}^{t}{c(x_\tau,u_\tau,\xi_\tau)}}{\ell}\right\}&\text{if} ~t+1 \in \mathcal{T'}, \\
\phi_{t}&\text{otherwise},
\end{cases}
    % \phi_{t+\ell}=\max\left(\phi_t,\frac{1}{\ell}\sum_{\tau=t}^{t+\ell-1}{h(\bm{x}_\tau,\bm{u}_\tau},\bm{\xi}_\tau)\right), \quad \forall t \in \mathbb{T}
\end{equation}
where $\mathcal{T'}=\{0,\ell,2\ell,\dots,T-\ell\}$ denotes the set of the beginnings of each measurement window. Note that $\phi_{t+1} \ge \phi_{t}$ for every $t \in \mathcal{T}$.
Therefore, the optimal DLS-DC \eqref{model:osocp} can be equivalently formulated as:
\begin{equation}
\label{model:eq_soc}    
  \begin{aligned}
\max_{\{\mu_t\}_{t \in \mathcal{T}}} ~& \mathbb{E} \left[ \sum_{t\in \mathcal{T}}{\mathcal{G}_t(\tilde{x}_t,\tilde{u}_t,\xi_t)} - \mathcal{C}(\tilde{\phi}_T) \right ] \\
\text{s.t.}~& \eqref{eqn:sys_state}\eqref{eqn:sys_cons}\eqref{eqn:sys_control}\eqref{eqn:sys_1}\eqref{eqn:new_state_update}, \\
%\mathbf{x}_{t+1} = \bm{f}_t(\bm{x}_t, \bm{\mu}_t, \bm{\xi}_t), \\
& \text{initializing}~\tilde{\phi}_0 = \phi_0,
% & \hspace{3.7cm} \mathbb{T}:=\{\ell,\dots,T-l,T\}. \nonumber
  \end{aligned}
\end{equation}
%\subsubsection{BMPC Formulation} % (fold)
%\label{ssub:block_mpc_formulation}
where $\tilde{\phi}_t$ denotes the decision variable corresponding to the new state at stage $t$. With the new state variable, it is much easier to apply the idea of MPC on the reformulated problem \eqref{model:eq_soc}. Generally, MPC considers the optimal control problem of a shorter control horizon $\{t,\cdots,t+W\}$ and utilizes a forecasted trajectory $\{\hat{\xi}_k\}_{k=t}^{t+W-1}$. As a result, BMPC solves the following optimal control problem for a rolling window of length $W$ at $t \in \mathcal{T'}$:
\begin{subequations}
\label{model:bmpc}
\begin{align}
\max_{\{\tilde{u}_k\}_{k=t}^{t+W-1}} ~& \sum_{k=t}^{t+W-1}  \mathcal{G}_k\big(\tilde{x}_k,\tilde{u}_k,\hat{\xi}_k\big) - \mathcal{H}(\tilde{\phi}_{t+W})  \\
\text{s.t.} ~& \tilde{x}_{k+1} = f_t(\tilde{x}_k,\tilde{u}_k, \hat{\xi}_k), \label{eqn:bmpc_sys}\\
& \tilde{\phi}_k~\text{updates as in}~\eqref{eqn:new_state_update}, \label{eqn:bmpc_update}\\
&h_t(\tilde{x}_k, \tilde{u}_k, \hat{\xi}_k) \leq 0,\label{eqn:bmpc_cons}\\
& \hspace{3cm} k=t,\dots,t+W-1,  \nonumber \\
& \text{initializing}~\tilde{x}_t=x_t, 
\tilde{\phi}_t=\phi_t, \label{eqn:bmpc_1}
% & k=t,t+\ell,\dots,t+W-\ell \label{eqn:bmpc_dc}
\end{align}
\end{subequations}
where $\mathcal{H}(\tilde{\phi}_{t+W})$ is the BMPC terminal cost (to be specified in Section \ref{sub:choosing_the_end_of_horizon_function_H}). The main difference between BMPC and the nominal MPC (NMPC) (see Section \ref{sub:mpc_approaches}) is the \emph{block} structure, illustrated in Fig. \ref{f:bmpc}. Instead of moving from stage $t$ to $t+1$, BMPC moves one block ($\ell$ stages) each time, \emph{i.e.}, from $t$ to $t+\ell$. The optimal controls of \eqref{model:bmpc} in the first block $\{\tilde{u}^*_t,\cdots,\tilde{u}^*_{t+\ell-1}\}$ will be implemented. Others $\{\tilde{u}^*_{t+\ell},\cdots,\tilde{u}^*_{t+W-1}\}$ are only advisory. 

The BMPC approach is summarized as Algorithm \ref{alg:block-MPC} below. Two factors affect the performance of BMPC: an initial guess on the maximum average consumption $\phi_0$, and a terminal cost $\mathcal{H}(\phi_{t+W})$.
An accurate estimate on $\phi_0$ can be obtained using external information, e.g., learning from historical data.
The choice of the terminal cost $\mathcal{H}(\phi_{t+W})$ lies at the heart of BMPC solution to DLS-DC.
%$\mathcal{H}(\phi_{t+W})$ heavily relies on the demand charge structure,
Detailed discussions and comparisons are in Section \ref{sub:choosing_the_end_of_horizon_function_H}. %More importantly, BMPC can be easily extended to any stochastic scheduling problems with demand charge (e.g. storage optimal control), not limited to systems of deferrable load.
%the optimal performance in some cases, numerical results show that the slightly deviated guess of $\phi_0$ leads to near-optimal solutions. \XG{any numerical results? Simulation results of one day should suffice.}.

\begin{algorithm}[htb]
  \begin{algorithmic}[1]
      \STATE \textbf{Initialization:} Initialize system with $x_0$, $\phi_0$, $W$ and $\ell$. 
    \FOR{$t \in \{0,\ell,\cdots,T-\ell\}$}
    \STATE Observe the system current state $x_t$, the up-to-date highest average consumption $\phi_t$ and the actual input $\xi_t$;
    \STATE Forecast the random inputs $\{\hat{\xi}_k\}_{k=t+1}^{t+W-1}$ over stage $t+1$ to $t+W-1$;
    \STATE Solve the BMPC Problem \eqref{model:bmpc};
    \begin{equation}
      \{u_k\}_{k=t}^{k+W-1} = \text{BMPC}(x_t,\phi_t,\{\hat{\xi}_k\}_{k=t}^{t+W-1})
    \end{equation}
    \STATE Take the the first $\ell$-stage controls from the solution to \eqref{model:bmpc}: 
    $\bm{U}_{t}^{*}=(u^{*}_t,\dots,u^{*}_{t+\ell-1})$;
    %\STATE Take $\bm{u}^{*}_t=\mu_t(\bm{x}_t,\bm{\xi}_{t})$ and update system vector to $\bm{x}_{t+1}$ according to \eqref{eqn:system}; 
      \FOR {$j \in \{t,t+1,\dots,t+\ell-1\}$}
      \STATE Observe system state $x_j$ and the actual input $\xi_j$;
      \STATE Take $u^{*}_j=\mu_j(x_j,\xi_j)$ and update system state to $x_{j+1}$ according to \eqref{eqn:system_state};
      \STATE Update $\phi_t$ according to \eqref{eqn:new_state_update}.\
      \ENDFOR
    \ENDFOR
  \end{algorithmic}
  \caption{BMPC under Demand Charge}
  \label{alg:block-MPC}
\end{algorithm}

\subsection{BMPC Terminal Cost $\mathcal{H}(\phi_{t+W})$}
\label{sub:choosing_the_end_of_horizon_function_H}
Intuitively, good choices of the terminal cost $\mathcal{H}(\phi_{t+W})$ should reflect the amortization of the demand charge in the current operating interval $t$. Some primitive forms of the terminal cost could be: %Some possible choices of $\mathcal{H}(\phi_{t+W})$ include
\begin{subequations}
\begin{align}
& \mathcal{H}(\phi_{t+W}) := \mathcal{C}(\phi_{t+W}),\quad t\in \mathcal{T'}, \label{eqn:EOH_demand_charge} \\
& \mathcal{H}(\phi_{t+W}) := \frac{W}{T} \mathcal{C}(\phi_{t+W}),\quad t\in \mathcal{T'}. \label{eqn:EOH_average_demand_charge}
\end{align}
\end{subequations}
However, these two choices perform poorly in practice because 
\eqref{eqn:EOH_demand_charge} imposes the demand charge \emph{over the entire control horizon $T$} on the rolling window. When $W \ll T$, the demand charge $\mathcal{C}(\phi_{t+W})$ would dominate the total stage reward of the rolling window. As a result, the solution in this setting will often be so conservative that schedulers would rather sacrifice most of the stage reward than incur a large demand charge cost. In addition, \eqref{eqn:EOH_demand_charge} fails to capture the fact that the demand charge is only posed for the peak consumption.

A slightly better choice is \eqref{eqn:EOH_average_demand_charge}, which splits the demand charge \emph{equally} among $T$ stages. This choice essentially assumes that the states of deferrable loads within each measurement window are almost identical, which is often not true in practice.

We propose a more judicious choice %The last equality holds true only in the case where the demand charge function $\mathcal{C}(\cdot)$ is linear, e.g., in the case of EV charging.
\begin{equation}
\label{eqn:bmpc_dc_general}
\mathcal{H}(\phi_{t+W}) := \mathcal{C}(\phi_{t+W}-\phi_{t}),~t \in \mathcal{T'}.
\end{equation}
The rationale behind \eqref{eqn:bmpc_dc_general} is twofold.
First, it is clear that $\phi_{t+W} \ge \phi_{t}$ always holds true according to \eqref{eqn:new_state_update}. If the peak consumption of the current rolling window $\{t,\cdots,t+W\}$ is no higher than the previous one ($\phi_{t+W} = \phi_{t}$), no additional cost should be considered, \emph{i.e.}, $\mathcal{H}(\phi_{t+W}) = 0$.
Additional cost occurs only when the peak consumption increases, \emph{i.e.}, $\phi_{t+W} > \phi_{t}$. 

The BMPC terminal cost can be further justified for power system applications where the demand charge cost is linear.  In this case, we have
% Second, we can decompose the demand charge cost $\mathcal{C}(\phi_{T})$ as
% , which happens at the end of the control horizon, 

% the summation of all terminal costs $\mathcal{H}(\phi_{t+W})$ resemble the demand charge function, i.e., 
% \begin{equation}
% \sum_{t \in \mathbb{T}} \mathcal{H}(\phi_{t+W}) = \sum_{t \in \mathbb{T}} \left( \mathcal{C}(\phi_{t+W}) - \mathcal{C}(\phi_{t}) \right) = \mathcal{C}(\phi_T) - \mathcal{C}(\phi_0)
% \end{equation}
% the excess amount of the previous performance peak should be considered. If the previous one still holds, the cost on demand charge should not be imposed at the current optimization. Since demand charge function $\mathcal{C}(\cdot)$ is commonly to be linear in power systems, we have the following equation:
%\begin{prop}
%\label{prop:demand_charge_decomposition}
%If the demand charge function $\mathcal{C}(\phi)$ is linear and non-decreasing with $\phi$, then 
\begin{equation}
\mathcal{C}(\phi_T) = \mathcal{C}(\phi_0) + \sum_{t \in \mathcal{T'}} \mathcal{H}(\phi_{t+W}) = \mathcal{C}(\phi_0) + \sum_{t \in \mathcal{T}} \mathcal{C}(\phi_{k+1}-\phi_k),
\end{equation}
% (2) For any demand charge function $\mathcal{C}(\phi)$,
% \begin{equation}
% \mathcal{C}(\phi_0) + \sum_{t \in \mathcal{T}} \mathcal{V}_k(\phi_{k+1}, \phi_k)  = \mathcal{C}(\phi_T).
% \end{equation}
%\end{prop}
which enables us to define a revised stage reward function that considers the cost of the demand charge at each stage $t \in \mathcal{T}$:
\begin{equation}
\label{eqn:new_stage_reward_general}
     \mathcal{V}_t(x_t,u_t,\xi_t,\phi_t):=
     \mathcal{G}_t(x_t,u_t,\xi_t)-\mathcal{C}(\phi_{t+1}-\phi_{t}).
\end{equation}
It is clear that \eqref{eqn:bmpc_dc_general} is a direct result of formulating BMPC using \eqref{eqn:new_stage_reward_general}. The equation above reveals that \eqref{eqn:bmpc_dc_general} embeds the demand charge cost, which occurs at the end of control horizon, into each stage as decomposed in \eqref{eqn:new_stage_reward_general}. Therefore, \eqref{eqn:bmpc_dc_general} effectively avoids the inferior performance by directly using the demand charge structure such as \eqref{eqn:EOH_demand_charge} or \eqref{eqn:EOH_average_demand_charge}.

\subsection{Related MPC Approaches}
\label{sub:mpc_approaches}
We summarize two related MPC approaches as benchmarks in our comparison studies. The first is the nominal MPC without demand charge, where no demand charge penalty is added to the objective. The second is a modification of EMPC \cite{risbeck2020economic} so that it applicable for general measurement window size $\ell$. %the economic MPC \cite{risbeck2020economic} that models demand charge explicitly but assumes that the measurement window matches the scheduling decision interval ($\ell=1$). Our description is an adaptation of EMPC in the context of DLS-DC in a stochastic setting.  

%\subsubsection{Nominal MPC (NMPC)} % (fold)
%\label{sub:nmpc}
%Instead of moving $\ell$ steps every time as in the BMPC approach, NMPC moves only one step at each time (see Fig. \ref{f:bmpc}), \emph{i.e.}, solving the following optimization problem consisting of $W$ stages at every stage $t \in \mathcal{T}$:
%\begin{equation}
%\begin{aligned}
%\max~ J_t^{\text{NMPC}} &:= \sum_{k=t}^{t+W-1}{\mathcal{G}(\tilde{x}_k,\tilde{u}_k,\hat{\xi}_k)} \\
%\text{s.t.}~& \eqref{eqn:bmpc_sys}\eqref{eqn:bmpc_cons},\\
%& \text{initializing}~\tilde{x}_t=x_t.
%\end{aligned}
%\end{equation}
%where $(\hat{D}_{i,k},\hat{T}_{i,k}^d)$ are the forecast states of a new EV arrived at charger $i$ at stage $t+k$. 
%Unlike BMPC, only the optimal control $\tilde{u}^*_t$ will be implemented. Others $\{\tilde{u}^*_{t+1},\cdots,\tilde{u}^*_{t+W-1}\}$ are only advisory.
%Note that NMPC does not take demand charge into consideration, which is another major difference from BMPC.

\subsubsection{Nominal MPC (NMPC) without Demand Charge}
\label{ssub:nmpc}
Instead of moving $\ell$ steps every time as in the BMPC approach, NMPC moves only one step at each time (see Fig. \ref{f:bmpc}), \emph{i.e.}, solving the following optimization problem consisting of $W$ stages at every stage $t \in \mathcal{T}$:
\begin{equation}
\begin{aligned}
\max~ J_t^{\text{NMPC}} &:= \sum_{k=t}^{t+W-1}{\mathcal{G}(\tilde{x}_k,\tilde{u}_k,\hat{\xi}_k)} \\
\text{s.t.}~& \eqref{eqn:bmpc_sys}\eqref{eqn:bmpc_cons},\\
& \text{initializing}~\tilde{x}_t=x_t.
\end{aligned}
\end{equation}
%where $(\hat{D}_{i,k},\hat{T}_{i,k}^d)$ are the forecast states of a new EV arrived at charger $i$ at stage $t+k$.
Unlike BMPC, only the optimal control $\tilde{u}^*_t$ will be implemented. Others $\{\tilde{u}^*_{t+1},\cdots,\tilde{u}^*_{t+W-1}\}$ are only advisory.
Note that NMPC does not take the demand charge into consideration, which is another major difference from BMPC.

%NMPC-$\alpha$ imposes a fraction of the estimated demand charge in each rolling-window optimization. Here we generalize \eqref{eqn:EOH_demand_charge}-\eqref{eqn:EOH_average_demand_charge} when $l=1$ and consider the following problem:
%Different from NMPC, MPC-TW appends a terminal cost as the weight of demand charge for each scheduling interval. Here we consider the form of \eqref{eqn:EOH_average_demand_charge} introduced in Section \ref{sub:choosing_the_end_of_horizon_function_H} and the case when $\ell=1$. At each stage $t \in \mathcal{T} $, MPC-TW solves the following problem:
%\begin{equation}
%\begin{aligned}
%\max~ J_t^{\text{MPC-}\alpha} &:= \sum_{k=t}^{t+W-1}{\mathcal{G}(\tilde{x}_k,\tilde{u}_k,\hat{\xi}_k)}-(1-e^{-\alpha t})\mathcal{C}(\phi_{t+W}) \\
%\text{s.t.}~& \eqref{eqn:bmpc_sys}\eqref{eqn:bmpc_update}\eqref{eqn:bmpc_cons},\\
%& \text{initializing}~\tilde{x}_t=x_t,\tilde{\phi}_t= \phi_t. 
%\end{aligned}
%\end{equation}
%Note that MPC-$\sigma$ reduces to NMPC when $\sigma_t=0$ for all $t \in \mathcal{T}$, and %it reduces to the approaches considered in \cite{kumar2018stochastic} with $\sigma_t=1$ and $W/T$ respectively.
%Note that the temporal structure of MPC-$\alpha$ is identical to NMPC as shown in Fig. \ref{f:bmpc}, and MPC-$\alpha$ reduces to NMPC when $\alpha=0$.  

\subsubsection{Economic MPC (EMPC) \cite{risbeck2020economic}}
\label{ssub:empc}
Here we demonstrate the formulation of EMPC under DLS-DC when $\ell=1$ and then present a slight modification of EMPC %demonstrate the formulation of EMPC under DLS-DC when $\ell=1$ and then 
so that it applies to cases when $\ell \geq 2$. %in \cite{Yang&Geng&Guan&Tong} . 

Let $(\mathbf{x}^{\text{ref}},\mathbf{u}^{\text{ref}},\bm{\xi}^{\text{ref}})$ be an arbitrarily known reference trajectory over $\mathcal{T}$. For each $t \in \mathcal{T}$, the objective of EMPC is defined as
\begin{equation}
\label{eqn:empc_obj}
  \begin{aligned}
   J^{\text{EMPC}}_t := \sum_{k=t}^{t+W-1}& {\mathcal{G}_k(x_k, u_k,\xi_k^{\text{ref}})} \\
  &-\mathcal{C}(\max(\phi_{t+W},\breve{\psi}_{t+W}^{\text{ref}})) - \mathcal{C}(\psi^{\text{ref}}),
  \end{aligned}
\end{equation}
where parameters $\psi^{\text{ref}}$ and $\breve{\psi}_{t+W}^{\text{ref}}$ denote the peak consumption over the entire and the remaining horizon (from stage $t+W$ to $T$) of the reference trajectory respectively, and are computed as 
%\begin{align}
%\label{eqn:empc_peak_new}
%&\psi^{\text{ref}} := \max_{t \in \mathcal{T'}}~\left\{ \frac{1}{\ell} \sum_{\tau = t }^{t+\ell-1}{c(x^{\text{ref}}_\tau,u^{\text{ref}}_\tau,\xi^\text{ref}_\tau)}  \right \}, \label{eqn:empc_peak_new} \\
%\label{eqn:empc_peaks_new}
%\breve{\psi}^{\text{ref}}_t &:=   
%\max_{k \in \mathcal{T'}, k \geq t}~\left\{\frac{1}{\ell}\sum_{\tau=k}^{k+\ell-1}{c(x_\tau^{\text{ref}},u_\tau^{\text{ref}},\xi^{\text{ref}}_\tau)}\right\}, ~t \in \mathcal{T}, \label{eqn:empc_peaks_new}
%\end{align}
\begin{align}
%\begin{equation}
%\label{eqn:empc_peak}
\psi^{\text{ref}} &:=\max_{t \in \mathcal{T}}~{c(x_t^{\text{ref}},u_t^{\text{ref}},\xi^{\text{ref}}_t)},  \label{eqn:empc_peak} \\
%\label{eqn:empc_peaks}
\breve{\psi}^{\text{ref}}_t &:=\max_{k \in \mathcal{T}, k \geq t}~{c(x_k^{\text{ref}},u_k^{\text{ref}},\xi^{\text{ref}}_k)}, ~t \in \mathcal{T}. \label{eqn:empc_peaks}
%\end{equation}
\end{align}

At stage $t$, EMPC can be formulated as
\begin{subequations}
\label{model:empc}
\begin{align}
\max_{\{\tilde{u}_k\}_{k=t}^{t+W-1}} ~& J_t^{\text{EMPC}} \nonumber \\
\text{s.t.}~& \tilde{x}_{k+1} = f_t(\tilde{x}_k, \tilde{u}_k,\xi_k^{\text{ref}}),  \label{eqn:empc_ref_state}  \\
& h_t(\tilde{x}_k,\tilde{u}_k,\xi_k^{\text{ref}}) \leq 0,   \label{eqn:empc_cons}\\
& \tilde{\phi}_{k+1}=\max(\tilde{\phi}_k,c(\tilde{x}_k,\tilde{u}_k,\xi_k^{\text{ref}})),  \label{eqn:empc_new_state} \\
&\hspace{2.5cm} k=t,\dots,t+W-1, \nonumber \\
& \tilde{x}_{t+W}=x_{t+W}^\text{ref}, \label{eqn:empc_state_cons} \\
& \text{initializing}~\tilde{x}_t=x_t,
\tilde{\phi}_t= \phi_t, \label{eqn:empc_1} 
%& \hspace{3.7cm} \mathcal{T'}:=\{0,\ell,2\ell,\dots,T-\ell\}. \nonumber
\end{align}
\end{subequations}
where \eqref{eqn:empc_new_state} represents a special case of \eqref{eqn:new_state_update} ($\ell=1$) and \eqref{eqn:empc_state_cons} the terminal constraint. The temporal structure of EMPC is the same as NMPC (see Fig. \ref{f:bmpc}).

%Now we show the implementation of EMPC when $\ell \geq 2$. 
%\paragraph{$\ell=1$}
%As described in Section \ref{sub:mpc_approaches}, EMPC requests additional parameters $\psi^{\text{ref}}$ and $\breve{\psi}^{\text{ref}}_t$ from the reference trajectory, which are defined as:
%\begin{align}
%\begin{equation}
%\label{eqn:empc_peak}
%\psi^{\text{ref}} &:=\max_{t \in \mathcal{T}}~{c(x_t^{\text{ref}},u_t^{\text{ref}},\overline{\xi})},  \label{eqn:empc_peak} \\
%\label{eqn:empc_peaks}
%\breve{\psi}^{\text{ref}}_t &:=\max_{k \in \mathcal{T}, k \geq t}~{c(x_k^{\text{ref}},u_k^{\text{ref}},\overline{\xi})}, ~t \in \mathcal{T}. \label{eqn:empc_peaks}
%\end{equation}
%\end{align}
%After computing the values of these parameters by \eqref{eqn:empc_peak} and \eqref{eqn:empc_peaks}, EMPC can be formulated as the standard form \eqref{model:empc}.
%\paragraph{$\ell \geq 2$} 
Although the original EMPC does not consider the case when the measurement window mismatches the scheduling interval ($\ell \geq 2)$, we only need to amend the values of the requested parameters from the reference trajectory so that EMPC can still be implemented. Since the demand charge is assessed according to the average consumption over $\ell$ consecutive stages, we re-compute parameters $\psi^{\text{ref}}$ and $\breve{\psi}^{\text{ref}}_t$ as 
\begin{align}
%\label{eqn:empc_peak_new}
&\psi^{\text{ref}} := \max_{t \in \mathcal{T'}}~\left\{ \frac{1}{\ell} \sum_{\tau = t }^{t+\ell-1}{c(x^{\text{ref}}_\tau,u^{\text{ref}}_\tau,\xi^\text{ref}_\tau)}  \right \}, \label{eqn:empc_peak_new} \\
%\label{eqn:empc_peaks_new}
\breve{\psi}^{\text{ref}}_t &:=   
\max_{k \in \mathcal{T'}, k \geq t}~\left\{\frac{1}{\ell}\sum_{\tau=k}^{k+\ell-1}{c(x_\tau^{\text{ref}},u_\tau^{\text{ref}},\xi^{\text{ref}}_\tau)}\right\}, ~t \in \mathcal{T}, \label{eqn:empc_peaks_new}
\end{align}
%where \eqref{eqn:empc_peak_new} and \eqref{eqn:empc_peaks_new} are the average peak consumption over each measurement window along the reference trajectory. 
Then we use these amended values for the input parameters as requested by EMPC.
%Two additional parameters from the reference trajectory are also required for the EMPC:
%\label{eqn:empc_peak}
%\psi^{\text{req}}=\max_{t \in \mathcal{T}}~{h(x_t^{\text{ref}},u_t^{\text{ref}},\overline{\xi})}
%\end{equation}
%\begin{equation}
%\label{eqn:empc_peaks}
%\breve{\psi}^{\text{req}}_t=\max_{k \in \mathcal{T}, k \geq t}~{h(x_k^{\text{ref}},u_k^{\text{ref}},\overline{\xi})}, ~t \in \mathcal{T}
%\end{equation}
%where \eqref{eqn:empc_peak} and \eqref{eqn:empc_peaks} compute the peak consumption of the entire and the remaining section of the reference trajectory respectively. Then, EMPC can be formulated as a standard form shown in \cite{risbeck2019economic}.
%\begin{equation}
%%\psi^{\text{new}}= \max_{t \in \mathcal{T'}}~\left\{ \frac{1}{\ell} \sum_{\tau = t }^{t+\ell-1}{h(x^{\text{ref}}_\tau,u^{\text{ref}}_\tau,\overline{\xi})}  \right \}
%\end{equation}
%\begin{equation}
%\label{eqn:empc_peaks_new}
%\breve{\psi}^{\text{new}}_t=   
%\max_{k \in \mathcal{T'}, k \geq t}~\left\{\frac{1}{\ell}\sum_{\tau=k}^{k+\ell-1}{h(x_\tau^{\text{ref}},u_\tau^{\text{ref}},\overline{\xi})}\right\}, ~t \in \mathcal{T}
%\end{equation}
%where \eqref{eqn:empc_peak_new} and \eqref{eqn:empc_peaks_new}  average the peak consumptions of the reference trajectory. 
%Then EMPC uses these new values as input parameters. 
However, it is clear that the recomputed values would not match those required by EMPC, which would cause a mismatch to the reference trajectory. %(see Section \ref{sub:5m}).
%Hence we keep the formulation unchanged and simply use the average power consumption instead of the stagewise, EMPC is formulated as 
%\begin{subequations}
%  \begin{align}
%   \max_{\{u_k\}_{k=t}^{t+W-1}} ~& \sum_{k=t}^{t+W-1}  \mathcal{G}_k\big(\tilde{x}_k,u_k,\hat{\xi}_k\big) - \mathcal{H}(\phi_{t+W})  \\
%\text{s.t.} ~& \tilde{x}_{k+1} = f(\tilde{x}_k,u_k, \hat{\xi}_k), \label{eqn:bmpc_sys}\\
%& \phi_k~\text{updates as in}~\eqref{eqn:new_state_update},\\
%&\omega(\tilde{x}_k, u_k, \hat{\xi}_k) \leq 0,\label{eqn:bmpc_cons}\\
%& \hspace{3cm} k=t,\dots,t+W-1  \nonumber \\
%& \tilde{x}_t=x_t \label{eqn:bmpc_1}
%  \end{align}
%\end{subequations}

BMPC differs from EMPC in the following aspects: 
\setcounter{paragraph}{0}
\paragraph{Scheduling Interval}
EMPC assumes that the resolution of the demand charge measurement matches the scheduling interval, while BMPC optimizes the demand charge according to the average consumption within multiple scheduling intervals, which fits in the practical cases.
\paragraph{Terminal Constraint and Terminal Cost} 
With known reference state and exogenous parameter trajectories, EMPC solves a deterministic multi-interval tracking problem under a terminal state constraint. However, EMPC does not accommodate cases when the exogenous parameter cannot be perfectly forecasted. %because it would deviate from the expected so that the hard state constraint becomes infeasible.
%to obtain strong performance guarantees, such a terminal constraint would be infeasible if considering the exogenous randomness. %However, the BMPC does not include such hard constraints thus there is no related feasibility issue resulting from this. 
EMPC also imposes a terminal cost to account for the demand charge, which needs nearly full information from a reference trajectory. %Since the average cost of the EMPC converges to a reference trajectory, the performance of this method is highly dependent on that reference. %However, the reference trajectory is always pre-determined before scheduling, hence it would not be practical to implement the method if the system states are influenced by exogenous stochastic variables. Whereas the Block MPC only needs an input parameter $\phi_0$, which can be approximated by an external problem as mentioned in Section \ref{alg:block-MPC}. 
In the absence of such a reference, it is unclear how to adjust the terminal cost to account accurately for demand charge. It should be noted that a simple adaptation of the terminal cost to $\mathcal{C}(\phi_{t+W})$ over penalizes the stage decisions. On the other hand, BMPC does not need to follow a reference trajectory, but to solve a small-scale deterministic optimization problem by exploiting generically short-term forecasts. Consequently, their terminal costs are quite different. %over its optimizations. This would lead to the domination of demand charge so that the economic performance would be significantly degraded.

\section{EV Charging via BMPC} % (fold)
\label{sec:application_ev_charging_scheduling}
We now specialize DLS-DC to tackle the problem of centralized scheduling of EV charging at public facilities \cite{Chen&Tong:12SGC,xu2012scheduling,huang2015matching,Yu&Chen&Tong16CSEE,yu2018deadline}. To this end, deferrable loads are EVs with charging demands that arrive stochastically, each with a random amount of charging need and specified deadline for completion\cite{yu2018deadline}. Here we adopt a Markov decision process (MDP) model widely used for the EV charging problems \cite{yu2018deadline,xu2016dynamic,j2019}.
%We study the optimal scheduling of EVs with stochastic arrival times, parking times and demands in this section.
%\subsection{Problem Formulation} % (fold)
%\label{sub:ev_charging}

 % subject to various constraints including maximum charging rate, time-varying electricity price and a monthly demand charge. 

\subsection{Nominal Model Assumptions} % (fold)
\begin{figure}[tb]
  	\centering
    \includegraphics[width=\linewidth]{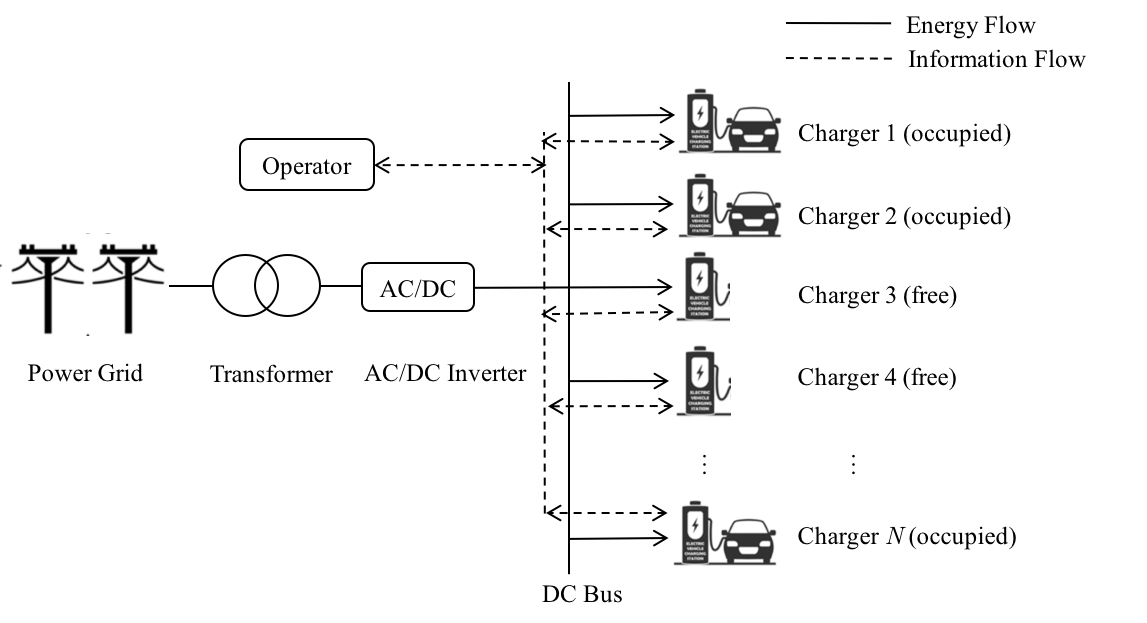}
	\caption{Schematic of EV charging scheduling}
 	\label{f:sys}
\end{figure}

\label{sub:nominal_model_assumptions}
% Some nominal model assumptions are presented below.
Consider an EV charging facility with $N$ chargers (charging ports) as illustrated in Fig. \ref{f:sys}. EVs arriving at the charging facility are assigned randomly to one of the available chargers. We assume that, upon arrival, the EV reveals to the operator its charging demand and deadline for completion. 
%\begin{figure}[htb]
%  	\centering
%    \includegraphics[width=\linewidth]{system.png}
%	\caption{perfect forecast}
% 	\label{f:sys}
%\end{figure}  

The operator faces a deadline scheduling problem, aimed at completing as many EV charging jobs as possible by their deadlines. The reward for the operator is the revenue from serving EV demands. The cost, on the other hand, comes from the electricity consumed in EV charging, the demand charge imposed by the distribution utility, and the penalty when the charging demand is not fulfilled. The operator also faces the constraint that only a finite number of chargers can be activated simultaneously due to transformer constraints from the distribution circuit.

Some of the key details of the EV charging problem and assumptions are outlined below.

\begin{description}
%\item[A1)] Time  (stage) is discrete, indexed by $t \in \mathcal{T} := \{0,1,\cdots,T-1\}$. Total $T$ stages $\{0,1,\cdots,T-1\}$ and a final stage $T$ are being considered.
  \item[A1)] All the chargers with a constant charging rate $R$ are available at any stage $t \in \mathcal{T}$. The charging decision at stage $t$ for the $i$th charger is a binary variable $u_{i,t} \in \{0,1\}$, with $1$ activating and $0$ deactivating the charging port. We also denote $M$ as the maximum number of simultaneous chargers allowed by the maximum power constraint of the local transformer, where $M<N$. 
  \item[A2)] The peak average consumption used to compute the demand charge is represented by variable $\psi$. In the case of EV charging, we assume the demand charge is a linear function $\mathcal{C}(\psi) = \pi^d \psi$ with the demand charge price $\pi^d$.
  \item[A3)] The EV arriving at the $i$th charger at the beginning of stage $t_0$ reveals random $D_{i,t_0}$ (the total amount of energy to be completed) and $T_{i,t_0}$ (the time for completion). %With probability $P(D,T^d)$, an EV with energy demand $D$ and parking time $T^d$ arrives at a available charger. With $P(0,0)$ the charger remains empty.  
  An EV will be automatically removed at the end of its completion time. An EV arriving at charger $i$ will be rejected if charger $i$ has been already occupied.
  \item[A4)] The operator receives a per unit reward $\pi^r$ and pays a time-varying charging cost $\pi^e_t$ if it serves an EV at stage $t$. For simplicity, we assume that parameters $\pi^r$ and $\{\pi^e_t\}_{t=0}^{T-1}$ are deterministic. The proposed approaches can be easily extended towards stochastic settings. 
  \item[A5)] If the total charging demand of EV $i$ is not completed at its completion time, then a penalty $qz_i$ occurs at price $q$, and $z_i$ denotes the amount of unmet demand. For simplicity, we assume that the penalty price is greater than the largest charging cost over the whole horizon, i.e., $q \ge \pi^e_t,~\forall t \in \mathcal{T}$.
 %%\item[A5)] The peak average consumption measured by demand charge is represented by variable $\psi$. In the case of EV charging, demand charge is a linear function $\mathcal{C}(\psi) = \pi^d \psi$ with demand charge price $\pi^d$.
\end{description}

\subsection{EV Charging as DLS-DC} % (fold)
\label{sub:EV_system}
We now define the DLS-DC model described in Section \ref{sec:DLS-DC} for the EV charging problem. 

\subsubsection{Exogenous Stochastic Input $\bm{\xi}=(\xi_{i,t})$}   
The input of DLS-DC model is a vector random process that models the arrivals of deferrable demands at individual chargers. The occupancy of each charger is an on-off process with the charger being occupied for the duration of the EV charging deadline and being idle for the duration of a Bernoulli process with parameter $p_i$ set by the overall arrival rate of the EV demand.  At the beginning of an occupied period of charger $i$, say at $t_0$, an EV arrives with random energy demand $D_{i,t_0}$ and random deadline $T_{i,t_0}$. Thus the input process at charger $i$ is given by $\xi_{i,t}= (D_{i,t_0},T_{i,t_0})$ for $t=t_0,\cdots, T_{i,t_0}$.  When the charger is idle, $\xi_{i,t}=(0,0)$. With probability $p_i$, $\xi_{i,t}=(0,0)$ transitions to $\xi_{i,t+1}=(D_{i,t+1}, T_{i,t+1})$.

\subsubsection{System State and State Evolution}
The state of charger $i$ at stage $t$ is given by a tuple $x_{i,t} = (r_{i,t} , \tau_{i,t})$, where $r_{i,t}$ represents the remaining demand to be served by deadline $T_{i,t}$ at charger $i$ and $\tau_{i,t}=T_{i,t}-t$ the lead time to the EV's completion at stage $t$. Hence, the system state is modeled as 
\begin{equation}
\label{eqn:EV_state}
x_{i,t+1} = \begin{cases} 
x_{i,t}-(u_{i,t},1) &~\text{if}~\tau_{i,t} > 1, \\
\xi_{i,t} & ~\text{if}~\tau_{i,t} \le 1.
\end{cases} 
\end{equation}
Note that when the charger is free, its state is $(0,0)$. When there is no EV arriving at charger $i$, the state of the charger remains at $(0,0)$.

\subsubsection{Constraints}
The total amount of power used for charging at one stage is limited by
\begin{equation}
\label{eqn:EV_power}
\sum_{i=1}^{N}{u_{i,t}} \le M,~ t \in \mathcal{T}.
\end{equation}
%Let $\mathcal{I}_t=\{i: t < t^d_i\}$ denote the set of occupied chargers at stage $t$, and $\mathcal{J}_t \subset \mathcal{I}_t$ is the set of EVs that will leave at stage $t+1$, i.e., $\mathcal{J}_t := \{i: T_i^d = t+1\} = \{i: \tau_{i,t}=1\}$. The total penalty $z_t$ incurred at stage $t$ is given by
%\begin{equation}
%\label{eqn:EV_penalty}
%z_t=\sum_{i \in \mathcal{J}_t}{(r_{i,t}-u_{i,t})}, \quad \forall t \in \mathcal{T}.
%\end{equation}
As shown in A2), the peak average demand $\psi$ over $\ell$  intervals is 
\begin{equation}
\label{eqn:EV_dc}
\psi=\max_{t \in \mathcal{T}'}
~{\frac{R}{\ell}\sum_{\tau = t}^{t+\ell-1}{\sum_{i=1}^{N}{u_{i,\tau}}}}. 
%\end{align}
\end{equation}

\subsubsection{Stage Reward}
The stage reward collected from all the EVs at stage $t$ is given by
\begin{equation}
\mathcal{G}_t(\mathbf{x}_t,\mathbf{u}_t,\bm{\xi}_t)= R\Delta\left[(\pi^r - \pi_t^e) \sum_{i=1}^{N} {u_{i,t}}-q\sum_{i \in \mathcal{J}_t}{(r_{i,t}-u_{i,t})}\right]  
\end{equation}
where $\Delta$ denotes the length of a scheduling stage and 
$\mathcal{J}_t$ the set of EVs that will leave at stage $t+1$, \emph{i.e.}, $\mathcal{J}_t : = \{i: \tau_{i,t}=1\}$.

\subsubsection{MDP formulation}
The objective of EV scheduling is to find the optimal control policy $\{\mu_t^*\}_{t \in \mathcal{T}}$ to maximize the expected total reward in the presence of the demand charge. At each stage $t$, a control law maps states to controls:
\begin{equation}
\label{eqn:policy_to_state}
\mathbf{u}_t = \mu_t(\mathbf{x}_t,\bm{\xi}_t).
\end{equation}
%where $\mathbf{x}_t=\{x_{1,t},\dots,x_{N,t}\}$ and $\mathbf{u}_t=\{u_{1,t},\dots,u_{N,t}\}$. $(\bm{D}_t,\bm{T^d}_t)$ represent the states of the newly arrived EVs at all chargers at stage $t$. 
Given an initial state $\mathbf{x}_0$, the EV schedule system can be formulated as  
\begin{equation}
\begin{aligned}
\label{model:EV}
\max_{\{\mu_t\}_{t \in \mathcal{T}}} ~& {\mathbb{E} \left[\sum_{t=0}^{T-1}{\mathcal{G}_t(\tilde{\mathbf{x}}_t,\tilde{\mathbf{u}}_t,\bm{\xi}_t)} -\pi^d\psi \right] } \\
\text{s.t.} ~& \eqref{eqn:EV_state}\eqref{eqn:EV_power}\eqref{eqn:EV_dc}\eqref{eqn:policy_to_state}, \\
& \text{initializing}~\tilde{\mathbf{x}}_0=\mathbf{x}_0.
\end{aligned}
\end{equation}
With \eqref{model:EV}, various MPC solutions, including the proposed BMPC approach, can be implemented.

\section{Numerical Results} % (fold)
\label{sec:numerical_results}
%We use the EV charging problem as a numerical example, and compare the performance of BMPC with four other methods (EMPC, NMPC, EDF, and LLF-LD). Sections \ref{sub:15m} and \ref{sub:5m} focus on the performance of BMPC and more discussions on other methods are in Section \ref{sec:discussions}.
%\subsection{Settings} % (fold)
%\label{sub:simulation_settings}
We conducted simulations involving stochastic EV-charging demands with random arrival times, charging demands and deadlines for completion.
We assumed that the number of newly arrived EVs at each stage followed a Poisson distribution. The charging demand and completion time of a new EV followed uniform distributions $\mathcal{U}(0,D_{\text{max}})$ and $\mathcal{U}(0,T_{\text{max}})$ respectively.

All prices were deterministic. The electricity prices were from the Electric Reliability Council of Texas (ERCOT)\footnote{Day-ahead Market (DAM) prices from November 1st to November 30th, 2019. Available at \url{http://www.ercot.com/mktinfo/dam}. }. The penalty price $q$ was set as 0.3 \$/kWh. In numerical simulations, we varied the demand charge price $\pi^d$ from 6 \$/kW to 21 \$/kW \cite{jin2016optimal}, where the length of the measurement window was fixed at 15 minutes and billing period a whole month. Other parameters are summarized in Table \ref{tab:parameters}.
\begin{table}[tb]
  \caption{Parameter settings} 
  \label{tab:parameters}
  \centering
  \begin{tabular}{lll}
    \toprule
    Parameter & Value & Note \\
    \midrule
     $N$ & 50 & Total Number of Chargers \\
     $R$   &  240 kW & Constant Charging Power\\
     $M$ & 25 & Maximum Number of Simultaneous Chargers \\
     $\lambda$ & 5 & Expected EV Arrival Rate \\
     $D_{\text{max}}$ & 120 kWh & Maximum Energy Demand \\
     $T_{\text{max}}^d$ & 1 hour & Maximum Completion Time \\
     %$W$ & 3 hours & Length of Rolling Window \\
    \bottomrule
  \end{tabular}
\end{table}
%Simulation results with the aforementioned settings are in Sections \ref{sub:15m} and \ref{sub:5m}. For a given demand charge price, we simulated all methods over 100 scenarios with randomly generated EV charging requests, and reported the average performance over these scenarios.
% subsection simulation_settings (end)

\subsection{Benchmark and Performance}
To compare with BMPC, we adopted both MPC-based approaches (NMPC and EMPC, see Section \ref{sub:mpc_approaches}) and index rules (EDF \cite{liu1973scheduling} and LLF-LD \cite{j2019}) as the benchmark methods. For each sampled trajectory, we ran each algorithm and computed its total reward gap to the upper bound (in percentage) as the performance measure.
%Although the charging demands and arrival/departure of EVs are stochastic in \eqref{model:osocp}, we can obtain an upper bound on each sampled trajectory of EV charging requests. 
Suppose we obtained $S$ sampled trajectories for all EV charging requests across time $\{\bm{\xi}_t^s,t=0,\dots,T\}_{s=1}^S$. We then solved an integer program that defined the \emph{deterministic} DLS-DC for the upper bound of the total reward on the $s$th trajectory:
\begin{equation}
\label{model:ub_each}
\begin{aligned}
\max_{\{\tilde{\mathbf{u}}_t\}_{t\in \mathcal{T}}}~& \sum_{t=0}^{T-1}{\mathcal{G}_t(\tilde{\mathbf{x}}_t,\tilde{\mathbf{u}}_t,\bm{\xi}_t^s)} -\pi^d\psi \\
\text{s.t.} ~& \eqref{eqn:EV_state}\eqref{eqn:EV_power}\eqref{eqn:EV_dc},  \\
%& x_{i,t+1} = \begin{cases} 
%x_{i,t}-(u_{i,t},1) &~\text{if}~\tau_{i,t} > 1 \\
%(D^s_{i,t},T^{d,s}_{i,t}) & ~\text{if}~\tau_{i,t} \le 1
%\end{cases} 
& \text{initializing}~\tilde{\mathbf{x}}_0=\mathbf{x}_0. 
\end{aligned}
\end{equation}
By solving \eqref{model:ub_each}, we also obtained the optimum of the maximum average consumption measured by the demand charge for each trajectory. 
%Simulation results with the aforementioned settings and performance measures were in Sections \ref{sub:15m} and \ref{sub:5m}. 
For a given demand charge price, we simulated all methods over $S=100$ scenarios with randomly generated EV charging requests, and reported the average performances over these scenarios.

\subsection{Multi-resolution DLS-DC} % (fold)
\label{sub:5m}
In practice, the resolution of control can be significantly finer than that of the demand-charge measurement. For example, the measurement window size can be 15 minutes whereas the EV charging decisions can be made at the one to five minute resolution, i.e. $\ell=3 \sim 15$. The results presented in this section are from simulations with $\ell=3$, \emph{i.e.}, $\Delta=5$ minutes. 
%In real world, the demand charge measurement window is often a fixed length such as 15 minutes, while the scheduling of EV charging can have a much finer resolution such as 5 minutes. As discussed in Section \ref{sub:difficulties_of_dealing_with_demand_charge}, this induces a mismatch in timescale, i.e., $\ell$ charging actions to consider within each demand charge measurement window. This section studies this setting, and we assume that the demand charge measurement window is 15 minutes, which contains three 5-minute EV scheduling stages ($\ell=3$).
% The demand charge is calculated using the average charging power over 15 minutes.
\begin{figure}[tb]
	\centering
	\begin{subfigure}[t]{\linewidth}
  	\centering
    \includegraphics[width=\linewidth]{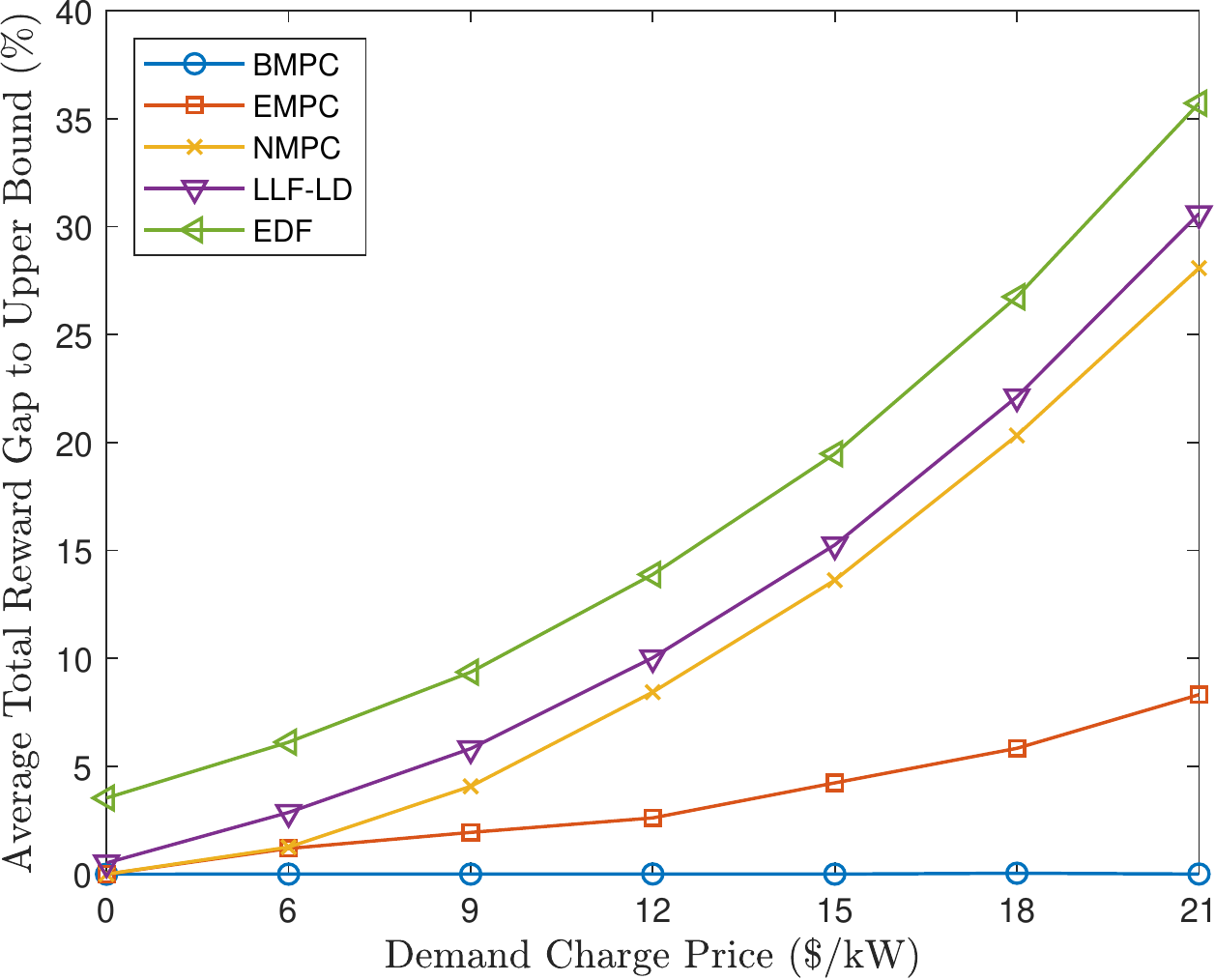}
	\caption{perfect forecast}
 	\label{f:5_det}
	\end{subfigure}  
    %\hspace{0.5cm}
	\begin{subfigure}[t]{\linewidth}
  	\centering
  	\vspace{0.5cm}
    \includegraphics[width=\linewidth]{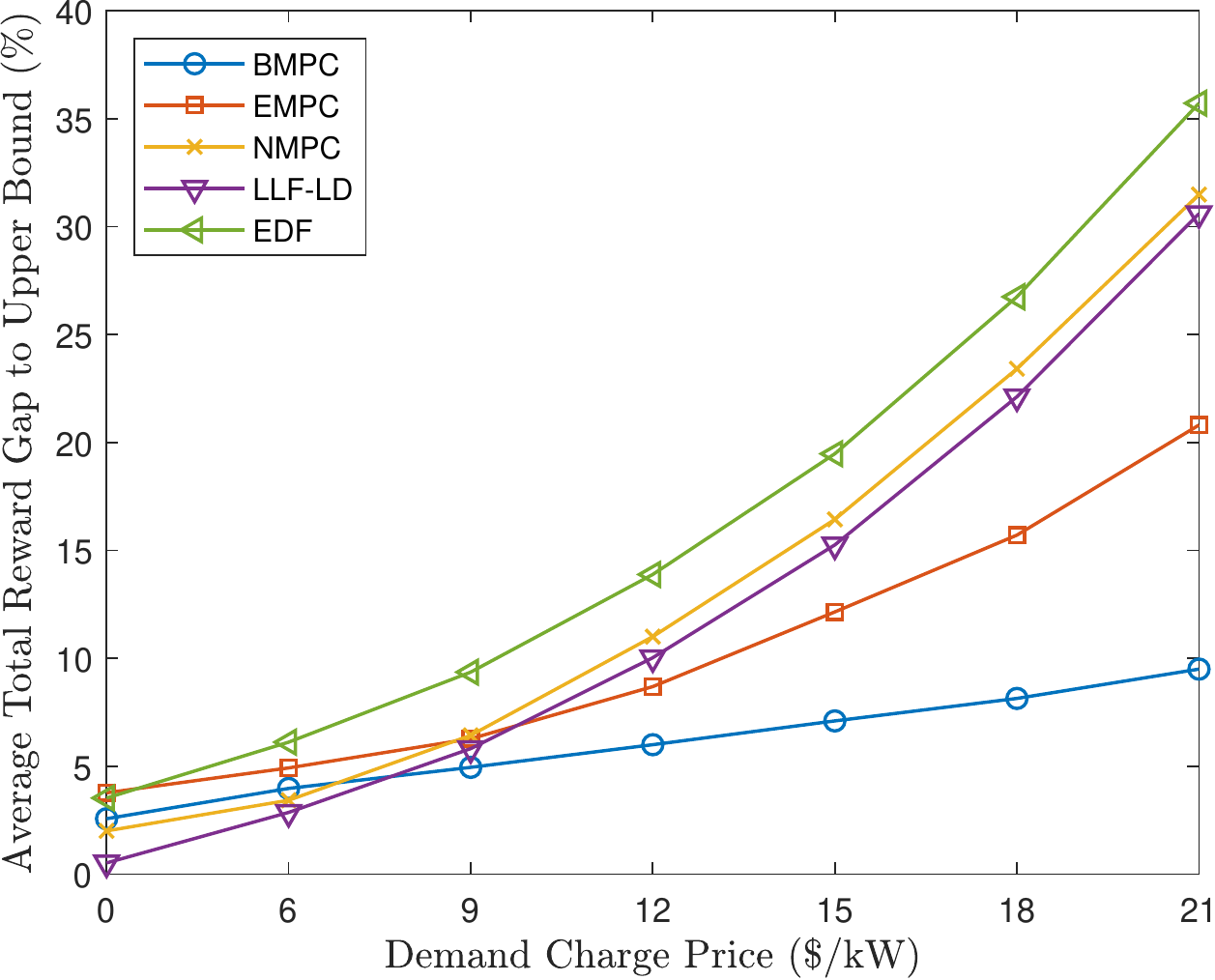}
	\caption{imperfect forecast}
 	\label{f:5_sto}
	\end{subfigure}  
	\caption{Average performance gap to upper bound with multi-resolution scheduling ($\ell=3$).}
	\label{fig:5min}
\end{figure} 

\subsubsection{Perfect Forecast} % (fold)
\label{ssub:perfect_forecast_5min}
% Similar to Section \ref{ssub:perfect_forecasts_15min},
In this case, we assumed that accurate information on EV arrivals, charging demands and required completion times were available, e.g., via reservation apps. 

Fig. \ref{f:5_det} compared the performance of BMPC with other methods under different demand charge prices (x-axis). The y-axis of Fig. \ref{f:5_det} showed the average optimality gaps between each method and the upper bound over all the sampled trajectories.  
%We first considered the case with accurate information on EVs. As shown in Fig. \ref{f:5_det}, 
We observed that BMPC outperformed other methods that did not consider the demand charge. In particular, BMPC almost reached upper bound ($0.16\%$ gaps), and achieved 27\% higher reward on average than NMPC. EMPC, however,
%Comparing Fig. \ref{f:5_det} with Fig. \ref{f:15_det}, we observed that EMPC 
performed worse than BMPC when $\ell = 3$. As mentioned in Section \ref{sub:mpc_approaches}, this was due to the mismatch between the actual peak values of the reference trajectory and those requested by EMPC. As shown in Table \ref{tab:power_5},  we observed that although both BMPC and EMPC managed to reduce the peak consumption as the demand charge price went higher, the peak consumption of EMPC deviated from the optimal one, which degraded its performance by at most 8\% compared to BMPC. 
% \begin{figure}[tb]
%   \centering
%     \includegraphics[width=\linewidth]{5_det.eps}
% 	\caption{Average performance to upper bound with 5-minute scheduling interval under perfect forecast}
%  \label{f:5_det}
% \end{figure}  
\begin{table}[tb]
  \caption{Maximum average power consumption (MW) with multi-resolution scheduling ($\ell=3$) under forecast (Others: NMPC, LLF-LD and EDF)} 
  \label{tab:power_5}
  \centering
  \begin{tabular}{ccccccccccc}
    \toprule
    \multicolumn{2}{c}{DC price (\$/MW)} & 0 & 6 &  9 & 12 & 15 & 18 & 21 \\
    \midrule
    %Upper Bound & 6.00 & 4.32 & 4.08 & 4.08 & 3.84 & 3.60 & 3.36 \\
    \multicolumn{2}{c}{Optimal} & 6.00 & 4.96 & 4.88 & 4.56 & 4.56 & 4.32 & 4.24 \\
    \midrule
    \parbox[t]{2mm}{\multirow{3}{*}{\rotatebox[origin=c]{90}{\scriptsize{perfect}}}}~\parbox[t]{2mm}{\multirow{3}{*}{\rotatebox[origin=c]{90}{\scriptsize{forecast}}}} & BMPC & 6.00 & 4.96 & 4.88 & 4.56 & 4.56 & 4.32 & 4.24\\
    & EMPC & 6.00 & 5.28 & 5.04 & 5.04 & 5.04 & 4.80 & 4.80\\
    & Others & 6.00 & 6.00 & 6.00 & 6.00 & 6.00 & 6.00 & 6.00\\  
    \midrule
    %Upper Bound & 6.00 & 4.32 & 4.08 & 4.08 & 3.84 & 3.60 & 3.36 \\
    \parbox[t]{2mm}{\multirow{3}{*}{\rotatebox[origin=c]{90}{\scriptsize{imperfect}}}}~\parbox[t]{2mm}{\multirow{3}{*}{\rotatebox[origin=c]{90}{\scriptsize{forecast}}}} & BMPC & 6.00 & 4.96 & 4.88 & 4.56 & 4.56 & 4.32 & 4.24\\
    & EMPC & 5.04 & 5.04 & 5.04 & 5.04 & 5.04 & 5.04 & 5.04\\
    & Others & 6.00 & 6.00 & 6.00 & 6.00 & 6.00 & 6.00 & 6.00\\  
    \bottomrule
  \end{tabular}
\end{table}
% \begin{table}[tb]
%   \caption{Average power peak (MW) with 5-minute scheduling interval under perfect forecast (Others: NMPC, LLF-LD and EDF)} 
%   \label{tab:power_5_det}
%   \centering
%   \begin{tabular}{p{1.6cm}cccccccc}
%     \toprule
%     DC price & 0 & 6 &  9 & 12 & 15 & 18 & 21 \\
%     \midrule
%     %Upper Bound & 6.00 & 4.32 & 4.08 & 4.08 & 3.84 & 3.60 & 3.36 \\
%     BMPC & 6.00 & 4.96 & 4.88 & 4.56 & 4.56 & 4.32 & 4.24\\
%     EMPC & 6.00 & 5.28 & 5.04 & 5.04 & 5.04 & 4.80 & 4.80\\
%     Others & 6.00 & 6.00 & 6.00 & 6.00 & 6.00 & 6.00 & 6.00\\  
%     \bottomrule
%   \end{tabular}
% \end{table}

\subsubsection{Imperfect Forecast} % (fold)
\label{ssub:imperfect_forecast_5min}
We considered a slightly more complicated setting, where the EV arrivals, charging demands and required completion times were random. For BMPC and NMPC, the forecasts on exogenous process $\bm{\xi}$ were based on the mean trajectories. The computation of forecasted reference trajectory of EMPC was demonstrated in Appendix \ref{sec:benchmarks}.%{Yang&Geng&Guan&Tong}. %\ref{sec:benchmarks}.
%We considered imperfect forecasts on EVs in this section.
%Predictions (EV arrivals, charging demands, and required completion times) were generated using the same method as shown in Section \ref{ssub:imperfect_forecasts_15m}. %The reference trajectories for EMPC were still obtained by solving \eqref{model:empc_ref}.

It was worth noting that BMPC may need to modify the schedule at certain stages, since it took a block of controls based on the inaccurate information for the near future. For example, BMPC would commit to the charging actions $\{\mathbf{\tilde{u}}^*_t,\mathbf{\tilde{u}}^*_{t+1},\cdots,\mathbf{\tilde{u}}^*_{t+\ell-1}\}$ after solving \eqref{model:bmpc} at stage $t$. The subsequent actions from $\mathbf{\tilde{u}}^*_t$, which were optimal for the predicted EV trajectory, might become infeasible for the realized EV profile. One simple solution to this issue was to deactivate the chargers that were conducting the infeasible actions whenever such rescheduling was necessary. 

Similar with Fig. \ref{f:5_det}, Fig. \ref{f:5_sto} quantified the average gaps between all methods and the upper bound. %Due to prediction errors, all MPC-based methods had positive gaps.
Due to the potentially suboptimal rescheduling actions of BMPC, LLF-LD achieved the best performance when the demand charge was small (e.g., $\pi^d \le 6$ \$/kW in Fig. \ref{f:5_sto}). In such a regime, the reward from scheduling played a more important role than the demand charge cost, and the peak differences among all the methods were relatively small (see Table \ref{tab:power_5}), which made the MPC-based methods less average total reward due to prediction errors. When the demand charge price was relatively high, the savings on the demand charge that BMPC achieved dominated the penalties due to the rescheduling, thus BMPC outperformed other methods and achieved nearly 20\% more average total reward than LLF-LD at 21 \$/kW. Meanwhile, the performance of EMPC further downgraded to at most 10\% gap to BMPC and 20\% to the upper bound, since both of mismatching in peak information and following an inaccurate reference trajectory came into effect.

\subsection{Single-resolution DLS-DC} % (fold)
\label{sub:15m}
This section validates the case when the resolution of peak-consumption measurement matched that of the decision, \emph{i.e.} $\ell=1$ ($\Delta=15$ minutes). In this case, BMPC operated at the same timescales as NMPC and EMPC.
\begin{figure}[tb]
	\centering
	\begin{subfigure}[t]{\linewidth}
  	\centering
    \includegraphics[width=\linewidth]{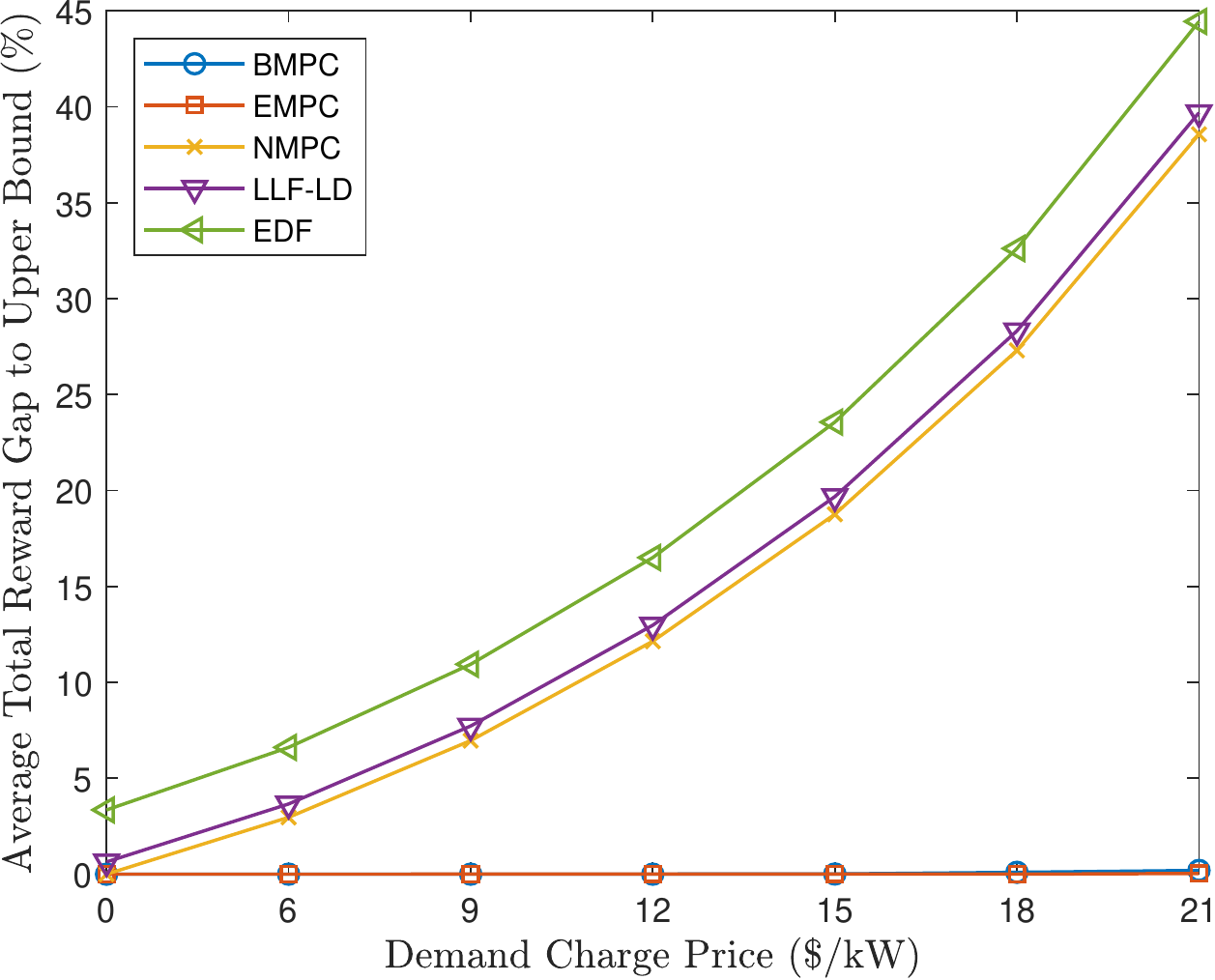}
	\caption{perfect forecast}
 	\label{f:15_det}
	\end{subfigure}
	%\hspace{0.5cm}
	\begin{subfigure}[t]{\linewidth}
  	\centering
  	\vspace{0.5cm}
    \includegraphics[width=\linewidth]{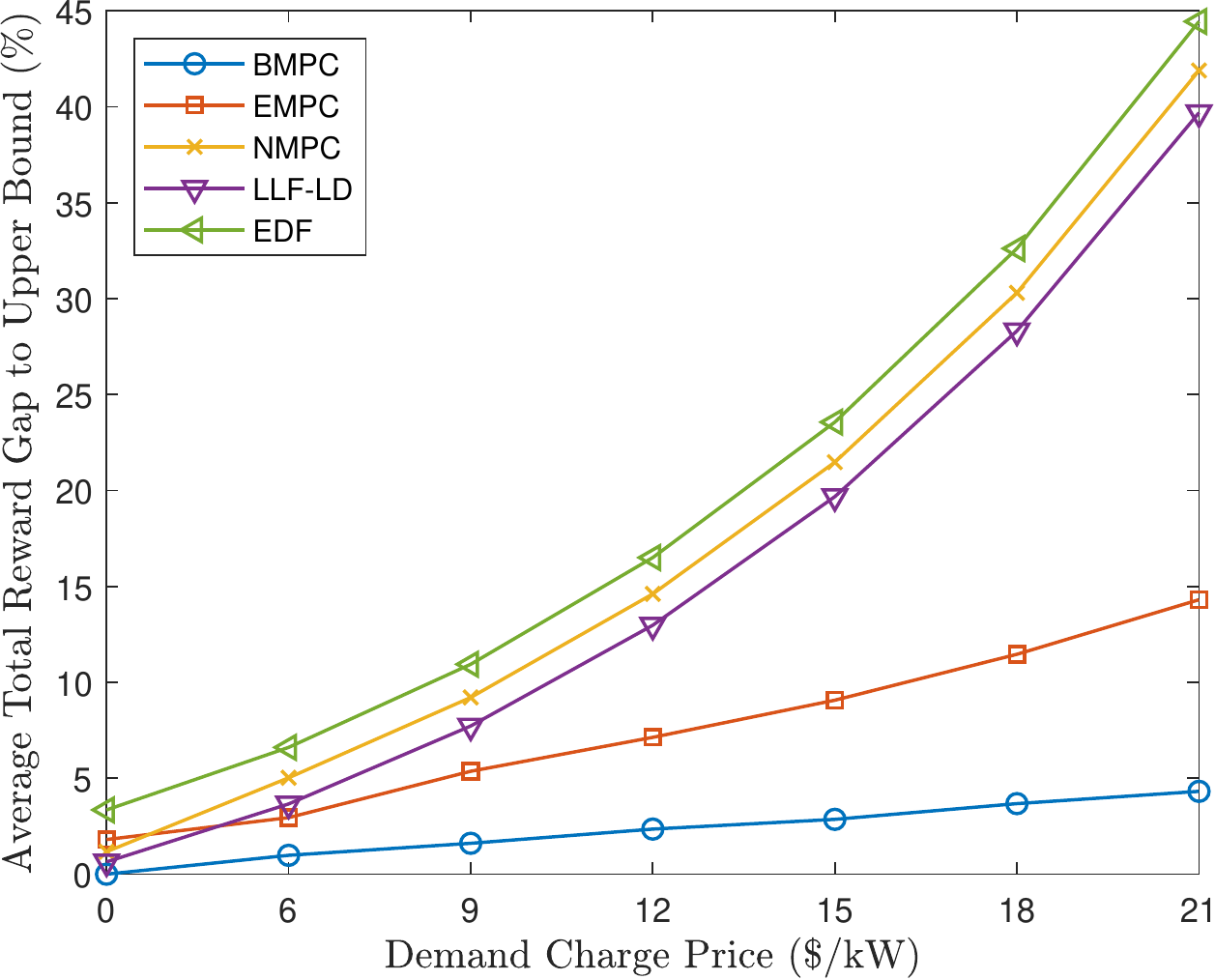}
	\caption{imperfect forecast}
 	\label{f:15_sto}
	\end{subfigure}  
	\caption{Average performance gap to upper bound with single-resolution scheduling ($\ell=1$).}
	\label{fig:15min}
\end{figure} 

\begin{table}[tb]
  \caption{Maximum average consumption (MW) with single-resolution scheduling ($\ell=1$) under forecast (Others: NMPC, LLF-LD and EDF)}
  \label{tab:power_15}
  \centering
  \begin{tabular}{ccccccccccc}
    \toprule
    \multicolumn{2}{c}{DC price (\$/MW)} & 0 & 6 &  9 & 12 & 15 & 18 & 21 \\
    \midrule
    \multicolumn{2}{c}{Optimal} & 6.00 & 4.32 & 4.08 & 4.08 & 3.84 & 3.60 & 3.36 \\
    \midrule
    %Upper Bound & 6.00 & 4.32 & 4.08 & 4.08 & 3.84 & 3.60 & 3.36 \\
    \parbox[t]{2mm}{\multirow{3}{*}{\rotatebox[origin=c]{90}{\scriptsize{perfect}}}}~\parbox[t]{2mm}{\multirow{3}{*}{\rotatebox[origin=c]{90}{\scriptsize{forecast}}}} & BMPC & 6.00 & 4.32 & 4.08 & 4.08 & 3.84 & 3.60 & 3.36\\
    & EMPC & 6.00 & 4.32 & 4.08 & 4.08 & 3.84 & 3.60 & 3.36\\
    %& MPC-$\alpha$ & 6.00 & 6.00 & 6.00 & 6.00 & 6.00 & 6.00 & 6.00 \\
    & Others & 6.00 & 6.00 & 6.00 & 6.00 & 6.00 & 6.00 & 6.00\\  
    \midrule
    %Upper Bound & 6.00 & 4.32 & 4.08 & 4.08 & 3.84 & 3.60 & 3.36 \\
    \parbox[t]{2mm}{\multirow{3}{*}{\rotatebox[origin=c]{90}{\scriptsize{imperfect}}}}~\parbox[t]{2mm}{\multirow{3}{*}{\rotatebox[origin=c]{90}{\scriptsize{forecast}}}} & BMPC & 6.00 & 4.32 & 4.08 & 4.08 & 3.84 & 3.60 & 3.36\\
    & EMPC & 5.04 & 4.32 & 4.32 & 4.32 & 4.32 & 4.32 & 4.32\\
    %& MPC-$\alpha$ & 6.00 & 6.00 & 6.00 & 6.00 & 6.00 & 6.00 & 6.00 \\
    & Others & 6.00 & 6.00 & 6.00 & 6.00 & 6.00 & 6.00 & 6.00\\  
    \bottomrule
  \end{tabular}
\end{table}

\subsubsection{Perfect Forecast} % (fold)
\label{ssub:perfect_forecasts_15min}

We first considered the case with accurate information on EVs.
%In this case, we assumed that accurate information on EV arrivals, charging demands and required completion times were available, e.g., via reservation apps.  
As shown in Fig. \ref{f:15_det},
%compared the performance of BMPC with other methods under different demand charge prices (x-axis). The y-axis of Fig. \ref{f:15_det} showed the average optimality gaps between each method and the upper bound over all the sampled trajectories. In this case, 
the best method was EMPC, which reached upper bound (almost $0\%$ gaps at all the demand charge prices) because it tracked the optimal reference trajectory. BMPC achieved similar performances but with slightly bigger gaps ($0.17\%$).
%As shown in Table \ref{tab:power_15}, the near-optimal performances of BMPC and EMPC were mainly due to the fact that they both managed to reduce the peak power over the month.
The other methods (NMPC, LLF-LD and EDF), which did not take the demand charge into account, all experienced rapid growth of the optimality gaps due to the large demand charge costs. %thus poor performances as LLF-LD and EDF in Fig. \ref{f:15_det}. When a considerable amount of total cost came from demand charge, e.g., $\pi^d=21$ \$/kW, LLF-LD gained $\sim37\%$ less reward than BMPC or EMPC. 

Table \ref{tab:power_15} further compared all methods in terms of peak consumption over the whole month. Both BMPC and EMPC reduced peak charging power as the demand charge prices increased, whereas the other methods failed to reduce the peak consumption thus reached the total charging limit 6 MW at all the demand charge prices.

% \begin{table}[tb]
%   \caption{Average power peak (MW) with 15-minute scheduling interval under perfect forecast} 
%   \label{tab:power_15_det}
%   \centering
%   \begin{tabular}{p{1.6cm}cccccccc}
%     \toprule
%     DC price & 0 & 6 &  9 & 12 & 15 & 18 & 21 \\
%     \midrule
%     %Upper Bound & 6.00 & 4.32 & 4.08 & 4.08 & 3.84 & 3.60 & 3.36 \\
%     BMPC & 6.00 & 4.32 & 4.08 & 4.08 & 3.84 & 3.60 & 3.36\\
%     EMPC & 6.00 & 4.32 & 4.08 & 4.08 & 3.84 & 3.60 & 3.36\\
%     NMPC, LLF-LD and EDF & 6.00 & 6.00 & 6.00 & 6.00 & 6.00 & 6.00 & 6.00\\  
%     \bottomrule
%   \end{tabular}
% \end{table}

\subsubsection{Imperfect Forecast} % (fold)
\label{ssub:imperfect_forecasts_15m}
We then considered imperfect forecasts on EVs.
Predictions were generated using the same method as shown in Section \ref{sub:5m}.
%We considered a slightly more complicated setting, where the EV arrivals, charging demands and required completion times were random. For BMPC and MPC-$\alpha$, the forecasts on exogenous process $\bm{\xi}$ are based on the mean trajectories. The computation of forecasted reference trajectory of EMPC was demonstrated in Appendix A.%\ref{sec:benchmarks}). %and we computed the actual total reward of EMPC decisions under realized trajectories.

%Similar with Fig. \ref{f:15_det}, Fig. \ref{f:15_sto} quantified the average gaps between all methods and the upper bound. 
Due to prediction errors, all MPC-based methods had positive gaps. 
We observed that BMPC was the best method with optimality gaps less than 5\%. EMPC, however, was no longer the best choice. EMPC achieved 10\% less reward than BMPC when the demand charge price was at 21 \$/kW. The main reason for the marked gap increase was that it tracked an untrustworthy reference trajectory. Since LLF-LD and EDF did not require predictions, their performances remained the same as the perfect forecast case.

\subsection{Discussions} % (fold)
\label{sub:other_methods}
Figs. \ref{fig:5min} and \ref{fig:15min} also illustrated the performance of NMPC, EDF, and LLF-LD. Comparing with BMPC and EMPC, the suboptimal performances of these methods were consequences of not considering the demand charge costs. 

As shown in all figures,  EDF was always the one with the least total reward in all cases. EDF was a simplistic and myopic scheduling policy; it failed to meet many charging requests thus suffered from large penalties. %Due to the small weight on demand charge, MPC-$\alpha$ failed to reduce demand charge as expected, and therefore had similar performances as NMPC and LLF-LD.

In the cases with perfect forecast (see Figs. \ref{f:5_det} and \ref{f:15_det}), NMPC performed slightly better than LLF-LD. In the cases with imperfect forecast (see Figs. \ref{f:5_sto} and \ref{f:15_sto}), LLF-LD surpassed NMPC. This was because predictions were not required to run LLF-LD.
It should be noted that, as an index-like policy, LLF-LD had much lower computation cost than MPC-based techniques, and it performed the best when the demand charge was relatively low. 
%It was also worth mentioning that LLF-LD was a much faster algorithm than MPC-TP and NMPC, and was a good candidate for online scheduling algorithm without considering demand charge.

%\XG{This paragraph is from your writing, why LLF-LD considers price variations?} Similarly to the 15-minute scheduling case, both LLF-LD and EDF tried to serve as many EVs as they can, while LLF-LD took into account the energy price variation and it generated less penalty than EDF, leading to a better average performance.
%This was the direct result of demand charge cost getting more dominant.

% subsection other_methods (end)

% section discussions (end)

\section{Conclusion} % (fold)
\label{sec:conclusion}
We consider the problem of DLS-DC, which can be widely adopted to applications such as scheduling of EV charging and cloud computing services. Due to the difficulties of multiple timescales posed by the demand charge pricing, we propose the BMPC algorithm with a special terminal cost to incorporate the demand charge at each scheduling stage. Through a motivating application of EV charging scheduling, our proposed approach shows advantageous performances compared to the benchmark methods, highlighting the significant impact of demand charge for deferrable load scheduling.

% section concluding_remarks (end)

\appendices
\section{Selected Benchmark Solutions}
\label{sec:benchmarks}

%\subsubsection{Nominal MPC (NMPC)} % (fold)
%\label{sub:nmpc}
%Instead of moving $\ell$ steps every time as in the BMPC approach, NMPC moves only one step at each time (see Fig. \ref{f:bmpc}), \emph{i.e.}, solving the following optimization problem consisting of $W$ stages at every stage $t \in \mathcal{T}$:
%\begin{equation}
%\begin{aligned}
%\max~ J_t^{\text{NMPC}} &:= \sum_{k=t}^{t+W-1}{\mathcal{G}(\tilde{x}_k,\tilde{u}_k,\hat{\xi}_k)} \\
%\text{s.t.}~& \eqref{eqn:bmpc_sys}\eqref{eqn:bmpc_cons},\\
%& \text{initializing}~\tilde{x}_t=x_t.
%\end{aligned}
%\end{equation}
%where $(\hat{D}_{i,k},\hat{T}_{i,k}^d)$ are the forecast states of a new EV arrived at charger $i$ at stage $t+k$.
%Unlike BMPC, only the optimal control $\tilde{u}^*_t$ will be implemented. Others $\{\tilde{u}^*_{t+1},\cdots,\tilde{u}^*_{t+W-1}\}$ are only advisory.
%Note that NMPC does not take demand charge into consideration, which is another major difference from BMPC.

\subsubsection{Reference Trajectory of Economic MPC \cite{risbeck2020economic}}
Here we introduce the computation of a reference trajectory for EMPC under stochastic settings. We assume the random inputs $\{\xi_t\}_{t=0}^{T-1}$ follow an arbitrary distribution with expectation at $\overline{\xi}$, then a forecasted reference trajectory can be computed by solving the following deterministic DLS-DC:
\begin{subequations}
\label{model:empc_ref}
\begin{align}
\max_{\{\tilde{u}_t\}_{t \in \mathcal{T}}} ~& \sum_{t\in \mathcal{T}}{\mathcal{G}_t(\tilde{x}_t, \tilde{u}_t,\overline{\xi})} - \mathcal{C}(\psi) \label{eqn:empc_ref_obj}\\
\text{s.t.}~& \tilde{x}_{t+1} = f_t(\tilde{x}_t, \tilde{u}_t,\overline{\xi}), \quad t \in \mathcal{T} \label{eqn:empc_ref_state_2}  \\
& h_t(\tilde{x}_t,\tilde{u}_t,\overline{\xi}) \leq 0, \quad t \in \mathcal{T}  \label{eqn:empc_ref_cons}\\
& u_t=\mu_t(\tilde{x}_t,\overline{\xi}), \quad t \in \mathcal{T} \label{eqn:empc_sys_control} \\
& \tilde{x}_0=x_0, \label{eqn:empc_ref_1} \\
&\psi = \max_{t \in \mathcal{T'}}~\left\{ \frac{1}{\ell} \sum_{\tau = t }^{t+\ell-1}{c(\tilde{x}_\tau,\tilde{u}_\tau,\overline{\xi})}  \right \}. \label{eqn:empc_ref_max_power_def} 
%& \hspace{3.7cm} \mathcal{T'}:=\{0,\ell,2\ell,\dots,T-\ell\}. \nonumber
\end{align}
\end{subequations}
By denoting $(\mathbf{x}^{\text{ref}},\mathbf{u}^{\text{ref}})$ as the optimal solution of \eqref{model:empc_ref}, the forecasted reference trajectory can be obtained as $(\mathbf{x}^{\text{ref}},\mathbf{u}^{\text{ref}},\overline{\bm{\xi}})$. %Now we introduce implementations of EMPC under two cases: 

\subsubsection{Earliest Deadline First (EDF)\cite{liu1973scheduling} } % (fold)
\label{ssub:earliest_deadline_first}
EDF is a rather simple online scheduling rule for deferrable loads. Specifically, at each stage $t$, it gives priorities to tasks with the earliest deadlines and tries to serve as many tasks as possible. Therefore, EDF would use the full power limit when the demand is heavy, resulting in large cost on demand charge.

\subsubsection{Least Laxity First with Later Deadline (LLF-LD) \cite{j2019}}% (fold)
\label{ssub:least_laxity_first_later_deadline_}
LLF-LD is an online algorithm for deferrable load scheduling, which prioritizes tasks with less laxity at each stage. Laxity, as defined in \cite{j2019}, is the difference between a server's lead time and its remaining processing time, reflecting the maximum number of stages that a task can tolerate before the time it has to be continuously processed to avoid non-completion penalty. For the EV charging problem in Section \ref{sec:application_ev_charging_scheduling}, the laxity of an EV at charger $i$ at stage $t$ is $\tau_{i,t}-r_{i,t}$. If the laxity of two tasks are the same, then it prioritizes the one with later deadline. Note that LLF-LD would also fully utilize the grid capacity as EDF, which results in significant amount of demand charge.

%{proof}[Proof of Proposition %\ref{prop:demand_charge_decomposition}]
%According to the definition of $\phi_t$ in \eqref{eqn:new_state_update}, we have
%\begin{equation}
%\phi_0 \le \phi_1 \le \cdots \le \phi_T 
%\end{equation}
%Since $\mathcal{C}(\phi_{k})$ are monotonically increasing in $\phi$:
%\begin{equation}
%\mathcal{C}(\phi_0) \le \mathcal{C}(\phi_1) \le \cdots \le \mathcal{C}(\phi_T)
%\end{equation}
%then $\mathcal{V}_k(\phi_{k+1}, \phi_k) = \mathcal{C}(\phi_{k+1}) - \mathcal{C}(\phi_{k})$, thus 
%\begin{multline*}
%\sum_{t \in \mathcal{T}} \mathcal{V}_k(\phi_{k+1}, \phi_k) = \sum_{k \in \mathcal{T}} \left( \mathcal{C}_k(\phi_{k+1}) - \mathcal{C}(\phi_k) \right) \\ = \mathcal{C}(\phi_T) - \mathcal{C}(\phi_0).
%\end{multline*}
%\end{proof}

\bibliographystyle{IEEEtran}
\bibliography{references}

\end{document}